\documentclass{amsart}

\usepackage{amsmath}
\usepackage[latin1]{inputenc}
\usepackage{amsfonts}
\usepackage{amssymb}
\usepackage{graphicx,epsfig}
\usepackage{amsthm}
\usepackage{amscd}
\usepackage[all]{xy}
\usepackage[latin1]{inputenc}

\usepackage[T1]{fontenc}

\usepackage{amsmath}

\usepackage{amsfonts}

\usepackage{amssymb}

\usepackage{amsthm}

\usepackage{graphics}

\newcommand{\mk}{\medskip}

\newcommand{\ZZ}{\mathbb{Z}}

\newcommand{\CC}{\mathbb{C}}

\newcommand{\NN}{\mathbb{N}}

\newcommand{\QQ}{\mathbb{Q}}

\newcommand{\Glie}{\mathfrak{g}}

\newcommand{\Hlie}{\mathfrak{h}}

\newcommand{\demo}{\noindent {\it \small Proof:}\quad}

\newcommand{\U}{\mathcal{U}}

\newcommand{\Lo}{\mathcal{L}}

\newtheorem{thm}{Theorem}[section]

\newtheorem{defi}[thm]{Definition}

\newtheorem{cor}[thm]{Corollary}

\newtheorem{prop}[thm]{Proposition}

\newtheorem{lem}[thm]{Lemma}

\usepackage[all]{xy}
\title{The Kirillov-Reshetikhin conjecture and solutions of $T$-systems}

\author{David Hernandez}
\address{\'Ecole Normale Sup\'erieure - DMA, 45, Rue d'Ulm F-75230 PARIS,
  Cedex 05  FRANCE}
\email{dhernand@dma.ens.fr URL: http://www.dma.ens.fr/$\sim$dhernand}

\begin{document}

\begin{abstract} We prove the Kirillov-Reshetikhin conjecture for all untwisted quantum affine algebras : we prove that the character
of Kirillov-Reshetikhin modules solve the $Q$-system and we give an explicit formula for the character of their tensor products. In the
proof we show that Kirillov-Reshetikhin
modules are special in the sense of monomials and that their $q$-characters solve the $T$-system (functional relations appearing in the study of solvable lattice models). Moreover we prove that the $T$-system can be written in the form of an exact sequence. 
For simply-laced cases, these results were proved by Nakajima in \cite{Nab, Nad} with geometric arguments (main result of \cite{Nab}) which are not available in general. The proof we use is different and purely algebraic, and so can be extended
uniformly to non simply-laced cases.

\vskip 4.5mm


\end{abstract}

\maketitle

\section{Introduction} Let $\U_q(\hat{\Glie})$ be an untwisted quantum affine algebra of rank $n$. The Kirillov-Reshetikhin 
modules form a certain infinite class of simple finite dimensional representations of $\U_q(\hat{\Glie})$. The main question answered in 
this paper is the following : what it is the character of the Kirillov-Reshetikhin modules and of their tensor products for the quantum group of finite type $\U_q(\Glie)\subset\U_q(\hat{\Glie})$ ? This problem goes back to 1931 as Bethe 
\cite{be} solved it for certain modules of type $A_1$ in another language. The methods to solve physical models involved here are now known as ``Bethe 
Ansatz''. In the 80s, in a serie of fundamental and striking papers, Kirillov and Reshetikhin \cite{kr, ki1, 
ki2, ki3} solved the problem for type $A$ and proposed formulas for all types by analyzing the Bethe Ansatz. These 
papers became the starting point of an intense research. The description of these characters is called the Kirillov-Reshetikhin conjecture. 

Let us describe the conjecture and recall some previous results on this problem in more details :

The Kirillov-Reshetikhin modules are simple $l$-highest weight modules (notion analog to the notion of highest weight 
module adapted to the Drinfeld realization of quantum affine algebras, see the definition \ref{lhigh}). They are 
characterized by their Drinfeld polynomials $(P_j(u))_{1\leq j\leq n}$ (analogs of the highest weight) which are of the form ($q_i=q^{r_i}$, see section \ref{un}):
$$P_i(u)=(1-au)(1-aq_i^2u)...(1-q_i^{2k-2}u)\text{ , and for $j\neq i$ : }P_j(u)=0$$ 
(for $k=1$ the Kirillov-Reshetikhin modules are called fundamental representations). The
Kirillov-Reshetikhin conjecture predicts the character of these modules 
and of their tensor products for 
the subalgebra $\U_q(\Glie)\subset\U_q(\hat{\Glie})$ ($\U_q(\Glie)$ is a quantum Kac-Moody algebra of finite type) : 
in \cite{kr} conjectural formulas were given for these characters (they were obtained by observation of the Bethe 
Ansatz 
related to solvable lattice models). A conjectural induction rule called $Q$-system (see the theorem \ref{conjdeux}) was also given in \cite{kr} (the 
$Q$-system for exceptional types was given in \cite{ki3}). We will denote by $\mathcal{F}(\nu)$ the formulas for the 
characters (see definition \ref{explicit}; we use here the version of \cite{ki1, ki2, hkoty, knt}. The version of \cite{kr} is 
slightly different because the definition of binomial coefficients is a little changed, see remark 1.3 of \cite{knt}). The Kirillov-Reshetikhin conjecture can be stated in the following form (see \cite{knt}) : the characters of Kirillov-Reshetikhin modules solve the $Q$-system and any of their tensors products are given by the formulas $\mathcal{F}(\nu)$. There are many results for these conjectures and related problems (see \cite{knt} for an historic and a guide through the huge literature on this subject; to name a few of very interesting recent ones see \cite{kl, hkoty, k, c0, Nab, Nad}), but no proof for all types.

\noindent In \cite{hkoty} it was proved that the $Q$-system implies the formulas $\mathcal{F}(\nu)$ (a certain asymptotic property, simplified in \cite{k}, is also needed). In \cite{kns} 
functional relations called $T$-system were defined (they are motivated by the observation of transfer matrices in 
solvable lattice models). Transposed in the language of Frenkel-Reshetikhin $q$-characters (analogs of characters for quantum affine 
algebras introduced in \cite{Fre}; see also \cite{kn}), and motivated by the relations between 
$q$-characters and Bethe-Anzatz in \cite{Fre}, it was naturally conjectured that the $q$-characters of 
Kirillov-Reshetikhin modules solve the $T$-system \cite{kosy} (see the theorem \ref{tsyst} for explicit definition of 
the $T$-systems). As in terms of usual characters the $T$-system becomes the $Q$-system, this conjecture combined with 
the results of \cite{hkoty, k} implies the Kirillov-Reshetikhin conjecture.

\noindent In general there is no explicit formula for $q$-characters of finite dimensional simple modules. However Frenkel and Mukhin \cite{Fre2} defined an algorithm to compute the $q$-character of a class of simple modules (satisfying the ``special'' property, term defined in \cite{Nab}, see below), and Nakajima \cite{Nab} defined an algorithm to compute the $q$-character of arbitrary simple modules in simply-laced cases. Although in general no explicit formula has been obtained from them (they are very complicated), in some cases they give useful informations.

\noindent In particular Nakajima \cite{Nab, Nad} made a remarkable advance by proving the Kirillov-Reshetikhin conjecture and $T$-systems in simply-laced cases : using the main result of \cite{Nab} (the algorithm) he proved that in simply-laced cases the Kirillov-Reshetikhin modules are special and noticed that this property is useful in the study of these modules (the algorithm of \cite{Nab} is drastically simplified in this situation). He also proved in \cite{Nad} that the $T$-systems that he established can be written in the form of exact sequences. 

\noindent The main result of \cite{Nab} is based on the study of quiver varieties \cite{Naams} and is not known in non simply-laced cases (see the conjecture of \cite{her02}), and so Nakajima's proof can not be used for all types.

In the present paper we propose a general uniform proof for all types of these conjectures (Kirillov-Reshetikhin conjecture, special property, $T$-systems and exact sequence). Our proof is purely algebraic (the results of \cite{Nab} are not used) and so can be extended to non simply-laced cases. 

Let us describe the results and methods of this paper in more details :

\noindent The $q$-characters morphism introduced by Frenkel and Reshetikhin \cite{Fre} for finite dimensional modules of quantum 
affine algebras is a ring morphism $\chi_q:\text{Rep}(\U_q(\hat{\Glie}))\rightarrow\ZZ[Y_{i,a}^{\pm}]_{i\in I, 
a\in\CC^*}$ that describes the decomposition of a representation in $l$-weight spaces for the Cartan subalgebra 
$\U_q(\hat{\Hlie})\subset\U_q(\hat{\Glie})$. So the monomials $m=\underset{i\in I, 
a\in\CC^*}{\prod}Y_{i,a}^{u_{i,a}(m)}$ of the Laurent polynomials ring $\ZZ[Y_{i,a}^{\pm}]_{i\in I, a\in\CC^*}$ are 
analogs of weight for quantum affine algebras. Note that we can get the usual character of a representation $V$ : it is 
equal to $(\beta\circ\chi_q)(V)$ where $\beta(m)=\underset{i\in I,a\in\CC^*}{\prod}e^{u_{i,a}(m)\Lambda_i}$ ($\Lambda_i$ is a 
fundamental weight). So to prove the Kirillov-Reshetikhin conjectures it suffices to get enough information on the 
$q$-character of Kirillov-Reshetikhin modules. Although one can not give a formula from it, the algorithm of Frenkel and Mukhin \cite{Fre2} works for special modules (ie. with a unique dominant monomial ($u_{i,a}(m)\geq 0$) in the corresponding $q$-character). So there are two points in our proof : first we prove that the Kirillov-Reshetikhin modules are special; then we extract enough information from this property (and the algorithm) to prove the $T$-system (and the Kirillov-Reshetikhin conjecture).

{\it First point : the Kirillov-Reshetikhin modules are special.} The problem of determining which modules are special is 
not easy in general, as we do not know a priori how to compute the $q$-character otherwise than having the explicit action 
of the Cartan subalgebra (the fundamental representations are special \cite{Fre, Fre2, Naams}, and in simply-laced cases the Kirillov-Reshetikhin modules are special \cite{Nab, Nad}). In \cite{her05} we proposed a new method to prove a certain ``cone'' property for simple representations of general quantum affinizations. This method is used in a fundamental way in this paper : it is based on an induction argument and the 
study of Weyl-module of type $sl_2$ (Weyl modules are some maximal finite dimensional representations of $l$-highest 
weight, see \cite{Cha4}). Adapting this argument, one gets information on what kind of monomials can appear in the 
$q$-character of a simple module, and prove that Kirillov-Reshetikhin modules are special for all types.

{\it Second point : the $q$-characters of Kirillov-Reshetikhin modules solve the $T$-system.} The "functional" 
relations of the $T$-system take place 
in $\text{Im}(\chi_q)$. This space is known to be $\underset{i\in I}{\bigcap}\text{Ker}(S_i)$ where the $S_i$ are the 
screening operators \cite{Fre, Fre2, Naams}. In particular an element of $\text{Im}(\chi_q)$ is characterized by his 
dominant monomials. So we have to know what are the dominant monomials in each member of the equality of the 
$T$-system, and what are 
the multiplicities. For this point we can get informations from the fact that the Kirillov-Reshetikhin modules are 
special : a product of Kirillov-Reshetikhin modules can only have some specific kinds of dominant monomials. With the 
help of an induction argument and the Frenkel-Mukhin algorithm, we achieve the goal of proving the conjectures.

{\it Complementary result : the exact sequence.} The $T$-system is a relation in the Grothendieck ring. As the category of 
finite dimensional representations of $\U_q(\hat{\Glie})$ is not semi-simple, it can not be directly transposed in 
terms of modules. However we get an exact sequence (generalizing \cite{Nad}) by proving that a certain tensor product of
Kirillov-Reshetikhin modules is 
simple, and by using a theorem of Chari \cite{c} and Kashiwara \cite{kas} that proves that a certain tensor product is 
$l$-highest weight.

For a geometric side, in analogy with the simply-laced cases, our result gives an explicit formula of what would be the
Euler number of a ``quiver variety'' in non simply-laced cases (a construction of such a variety is not known yet).

\noindent Note that the methods of our proof may also be used to prove the conjectures of \cite{hkott} for twisted cases, and to
establish a $T$-system for representations of general quantum affinizations (such as quantum toroidal algebras) considered in \cite{her04}. We let these points for another paper.

\noindent Note that for classical type $ABCD$, the formulas $\mathcal{F}(\nu)$ for the character of Kirillov-Reshetikhin modules can be written in another ``numerical'' form, see the theorem 7.1 of \cite{hkoty} (such formulas were first conjectured in \cite{kr}). So the results of this paper also imply these formulas.

Let us describe the organization of this paper :

In section \ref{un} we give backgrounds on quantum Kac-Moody algebras, Drinfeld realization, finite dimensional 
representations of quantum affine algebras and $q$-characters. In section \ref{deux} we give the Kirillov-Reshetikhin 
conjecture, $Q$-systems and $T$-systems : as proved in \cite{hkoty, k} the formulas $\mathcal{F}(\nu)$ (theorem \ref{conjun}) are a 
consequence of the $Q$-system (theorem \ref{conjdeux}). This $Q$-system is a consequence of the $T$-system (theorem 
\ref{tsyst}) that we establish in general by proving that the Kirillov-Reshetikhin modules are special (theorem \ref{domkr}) in section \ref{special}. The end of the proof of the theorem \ref{tsyst} is given in the 
section \ref{proofth}. In the section \ref{sequ} we prove that the $T$-system can be written in the form of an exact 
sequence (theorem \ref{exact}). In the section \ref{exptsyst} we give formulas of the $T$-system for each type.

{\bf Acknowledgments :} A part of this paper was written when the author visited the RIMS (Kyoto) in the summer of
2004. 
He is very grateful to Hiraku Nakajima for this invitation, for his help in the course of the preparation of this 
paper and for encouraging him to prove the Kirillov-Reshetikhin conjecture; his precious comments on a previous version of this paper improved the clarity 
of the proofs. The author would also like to thank Edward Frenkel, Nicolai Reshetikhin and Marc Rosso for their support.

\section{Background}\label{un}

\subsection{Cartan matrix and quantized Cartan matrix} Let $C=(C_{i,j})_{1\leq i,j\leq n}$ be a Cartan matrix of 
finite type. We denote 
$I=\{1,...,n\}$. $C$ is symmetrizable : there is a matrix $D=\text{diag}(r_1,...,r_n)$ ($r_i\in\NN^*$)\label{ri} such 
that $B=DC$\label{symcar} is symmetric. In particular if $C$ is symmetric then it $D=I_n$ (simply-laced case). 
\noindent We consider a realization $(\Hlie, \Pi, \Pi^{\vee})$ of $C$ (see \cite{kac}): $\Hlie$ is a $n$ dimensional $\QQ$-vector space, $\Pi=\{\alpha_1,...,\alpha_n\}\subset \Hlie^*$ (set of the simple roots) and $\Pi^{\vee}=\{\alpha_1^{\vee},...,\alpha_n^{\vee}\}\subset \Hlie$ (set of simple coroots) are set such that for $1\leq i,j\leq n$: $\alpha_j(\alpha_i^{\vee})=C_{i,j}$.
Let $\Lambda_1,...,\Lambda_n\in\Hlie^*$ (resp. $\Lambda_1^{\vee},...,\Lambda_n^{\vee}\in\Hlie$) be the the fundamental weights (resp. coweights) : $\Lambda_i(\alpha_j^{\vee})=\alpha_i(\Lambda_j^{\vee})=\delta_{i,j}$. Denote $P=\{\lambda \in\Hlie^*/\forall i\in I, \lambda(\alpha_i^{\vee})\in\ZZ\}$ the set of weights and $P^+=\{\lambda \in P/\forall i\in I, \lambda(\alpha_i^{\vee})\geq 0\}$ the set of dominant weights. For example we have $\alpha_1,...,\alpha_n\in P$ and $\Lambda_1,...,\Lambda_n\in P^+$. Denote $Q=\underset{i\in I}{\bigoplus}\ZZ \alpha_i\subset P$ the root lattice, $Q^+=\underset{i\in I}{\sum}\NN \alpha_i\subset Q$, $\Delta$ the set of roots and $\Delta^+$ the set of positive roots. For $\lambda,\mu\in \Hlie^*$, denote $\lambda \geq \mu$ if $\lambda-\mu\in Q^+$. Let $\nu:\Hlie^*\rightarrow \Hlie$ linear such that for all $i\in I$ we have $\nu(\alpha_i)=r_i\alpha_i^{\vee}$. For $\lambda,\mu\in\Hlie^*$, $\lambda(\nu(\mu))=\mu(\nu(\lambda))$.

\noindent In the following we suppose that $q\in\CC^*$ is not a root of unity. We denote $q_i=q^{r_i}$ and for $l\in\ZZ, r\geq 0, m\geq m'\geq 0$ we define in $\ZZ[q^{\pm}]$ :
$$[l]_q=\frac{q^l-q^{-l}}{q-q^{-1}}\in\ZZ[q^{\pm}]\text{ , }[r]_q!=[r]_q[r-1]_q...[1]_q\text{ ,
}\begin{bmatrix}m\\m'\end{bmatrix}_q=\frac{[m]_q!}{[m-m']_q![m']_q!}.$$

\noindent Let $C(z)$ be the quantized Cartan matrix defined by ($i\neq j\in I$): $$C_{i,i}(z)=z_i+z_i^{-1}\text{ ,
}C_{i,j}(z)=[C_{i,j}]_z$$ $C(z)$ is invertible (see \cite{Fre}). We denote by $\tilde{C}(z)$ the inverse matrix of
$C(z)$ and $D(z)$ the diagonal matrix such that for $i,j\in I$, $D_{i,j}(z)=\delta_{i,j}[r_i]_z$.

\subsection{Quantum algebras}\label{qkma}

\subsubsection{Quantum groups}

\begin{defi} The quantum group $\U_q(\Glie)$ is the $\CC$-algebra with generators $k_i^{\pm 1}$, $x_i^{\pm}$ ($i\in I$) and 
relations: 
$$k_ik_j=k_jk_i\text{ , } k_ix_j^{\pm}=q_i^{\pm C_{i,j}}x_j^{\pm}k_i,$$
$$[x_i^+,x_j^-]=\delta_{i,j}\frac{k_i-k_i^{-1}}{q_i-q_i^{-1}},$$
$$\underset{r=0... 1-C_{i,j}}{\sum}(-1)^r\begin{bmatrix}1-C_{i,j}\\r\end{bmatrix}_{q_i}(x_i^{\pm})^{1-C_{i,j}-r}x_j^{\pm}(x_i^{\pm})^r=0 \text{ (for $i\neq j$)}.$$
\end{defi}

\noindent This algebra was introduced independently by Drinfeld \cite{Dri1} and Jimbo \cite{jim}. It is remarkable that 
one can define a Hopf algebra structure on $\U_q(\Glie)$ by : 
$$\Delta(k_i)=k_i\otimes k_i,$$ 
$$\Delta(x_i^+)=x_i^+\otimes 1 + k_i\otimes x_i^+\text{ , }\Delta(x_i^-)=x_i^-\otimes 
k_i^{-1} + 1\otimes x_i^-,$$ 
$$S(k_i)=k_i^{-1}\text{ , }S(x_i^+)=-x_i^+k_i^{-1}\text{ , }S(x_i^-)=-k_ix_i^-,$$ 
$$\epsilon(k_i)=1\text{ , }\epsilon(x_i^+)=\epsilon(x_i^-)=0.$$

\noindent Let $\U_q(\Hlie)$ the commutative subalgebra of $\U_q(\Glie)$ generated by the $k_i^{\pm 1}$ ($i\in I$).

\noindent For $V$ a $\U_q(\Hlie)$-module and $\omega\in P$ we denote by $V_{\omega}$ the weight space of weight 
$\omega$: 
$$V_{\omega}=\{v\in V/\forall i\in I, k_i.v=q_i^{\omega(\alpha_i^{\vee})}v\}.$$ 
In particular we have $x_i^{\pm}.V_{\omega}\subset V_{\omega \pm \alpha_i}$.

\noindent We say that $V$ is $\U_q(\Hlie)$-diagonalizable if $V=\underset{\omega\in P}{\bigoplus}V_{\omega}$ (in particular $V$ is of type $1$).

\noindent For $V$ a finite dimensional $\U_q(\Hlie)$-diagonalizable module we set $\chi(V)=\underset{\omega\in P}{\sum}\text{dim}(V_{\omega})e^{\omega}\in\ZZ[e^{\omega}]_{\omega\in P}$ the usual character. 

\subsubsection{Quantum loop algebras} We will use the second realization (Drinfeld realization) of the quantum loop 
algebra $\U_q(\Lo\Glie)$ (subquotient of the quantum affine algebra $\U_q(\hat{\Glie})$) :

\begin{defi}\label{defiaffi} $\U_q(\Lo\Glie)$ is the algebra with
generators $x_{i,r}^{\pm}$ ($i\in I, r\in\ZZ$), $k_i^{\pm 1}$ ($i\in I$), $h_{i,m}$ ($i\in I, m\in\ZZ-\{0\}$) and the
following relations ($i,j\in I, r,r'\in\ZZ, m\in\ZZ-\{0\}$): 
$$\label{afcart}[k_i,k_j]=[k_{h},h_{j,m}]=[h_{i,m},h_{j,m'}]=0,$$
$$k_ix_{j,r}^{\pm}=q_i^{\pm C_{i,j}}x_{j,r}^{\pm}k_i,$$
$$[h_{i,m},x_{j,r}^{\pm}]=\pm \frac{1}{m}[mB_{i,j}]_qc^{|m|/2}x_{j,m+r}^{\pm},$$
$$[x_{i,r}^+,x_{j,r'}^-]=\delta_{i,j}\frac{\phi^+_{i,r+r'}-\phi^-_{i,r+r'}}{q_i-q_i^{-1}},$$
$$x_{i,r+1}^{\pm}x_{j,r'}^{\pm}-q^{\pm B_{i,j}}x_{j,r'}^{\pm}x_{i,r+1}^{\pm}=q^{\pm B_{i,j}}x_{i,r}^{\pm}x_{j,r'+1}^{\pm}-x_{j,r'+1}^{\pm}x_{i,r}^{\pm},$$
$$\underset{\pi\in\Sigma_s}{\sum}\underset{k=0..s}{\sum}(-1)^k\begin{bmatrix}s\\k\end{bmatrix}_{q_i}x_{i,r_{\pi(1)}}^{\pm}...x_{i,r_{\pi(k)}}^{\pm}x_{j,r'}^{\pm}x_{i,r_{\pi(k+1)}}^{\pm}...x_{i,r_{\pi(s)}}^{\pm}=0,$$
where the last relation holds for all $i\neq j$, $s=1-C_{i,j}$, all sequences of integers $r_1,...,r_s$. $\Sigma_s$ is
the symmetric group on $s$ letters. For $i\in I$ and $m\in\ZZ$, $\phi_{i,m}^{\pm}\in \U_q(\Lo\Glie)$ is determined
by the formal power series in $\U_q(\Lo\Glie)[[z]]$ (resp. in $\U_q(\Lo\Glie)[[z^{-1}]]$): 
$$\underset{m\geq 0}{\sum}\phi_{i,\pm m}^{\pm}z^{\pm m}=k_i^{\pm 1}\text{exp}(\pm(q-q^{-1})\underset{m'\geq
1}{\sum}h_{i,\pm m'}z^{\pm m'})$$ 
and $\phi_{i,m}^+=0$ for $m<0$, $\phi_{i,m}^-=0$ for $m>0$. \end{defi}

\noindent $\U_q(\Lo\Glie)$ has a structure of Hopf algebra (from the Hopf algebra structure of $\U_q(\hat{\Glie})$).

\noindent For $J\subset I$ we denote by $\Glie_J$ the semi-simple Lie algebra of Cartan matrix $(C_{i,j})_{i,j\in J}$. We 
have $\U_q(\Lo\Glie_J)\subset\U_q(\Lo\Glie)$ and for $i\in I$, $\U_q(\Lo\Glie_i)\simeq \U_{q_i}(\Lo sl_2)$.

\subsection{Finite dimensional representations of quantum loop algebras} Denote by 
$\text{Rep}(\U_q(\Lo\Glie))$ the Grothendieck ring of finite 
dimensional representations of $\U_q(\Lo\Glie)$.

\subsubsection{Monomials and $q$-characters}\label{defimono} The subalgebra $\U_q(\Lo\Hlie)\subset\U_q(\Lo\Glie)$ 
generated by the $\phi_{i,\pm 
m}^{\pm}, k_i^{\pm 1}$ is commutative, so for $V$ a $\U_q(\Lo\Glie)$-module we have :
$$V=\underset{\gamma=(\gamma_{i,\pm m}^{\pm})_{i\in I, m\geq 0}}{\bigoplus}V_{\gamma}$$
$$\text{where : }V_{\gamma}=\{v\in V/\exists p\geq 0, \forall i\in I, m\geq 0, (\phi_{i,\pm m}^{\pm}-\gamma_{i,\pm m}^{\pm})^p.v=0\}.$$
One can prove \cite{Fre} that $\gamma$ satisfies in $\CC[[u]]$ (resp. in $\CC[[u^{-1}]]$) :
\begin{equation}\label{formfin}\underset{m\geq 0}{\sum}\gamma_{i,\pm m}^{\pm}u^{\pm 
m}=q^{\text{deg}(Q_i)-\text{deg}(R_i)}\frac{Q_i(uq_i^{-1})R_i(uq_i)}{Q_i(uq_i)R_i(uq_i^{-1})}.\end{equation}
The Frenkel-Reshetikhin $q$-characters morphism $\chi_q$ \cite{Fre} encodes the $l$-weight $\gamma$ of $\U_q(\Lo\Hlie)$ (see also \cite{kn}). It is an injective ring morphism : 
$$\chi_q:\text{Rep}(\U_q(\Lo\Glie))\rightarrow \ZZ[Y_{i,a}^{\pm}]_{i\in I, a\in\CC^*}\text{ , }\chi_q(V)=\underset{\gamma}{\sum}\text{dim}(V_{\gamma})m_{\gamma}$$
where 
$$m_{\gamma}=\underset{i\in I, a\in\CC^*}{\prod}Y_{i,a}^{q_{i,a}-r_{i,a}}\text{ , }Q_i(u)=\underset{a\in\CC^*}{\prod}(1-ua)^{q_{i,a}}\text{ , }R_i(u)=\underset{a\in\CC^*}{\prod}(1-ua)^{r_{i,a}}.$$
The $m_{\gamma}$ are called monomials or $l$-weight (they are analogs of weight) and we denote $V_{\gamma}=V_{m_{\gamma}}$. We denote by $A$ the set of monomials of $\ZZ[Y_{i,a}^{\pm}]_{i\in I, a\in\CC^*}$.

\noindent For $J\subset I$, $\chi_q^J$ is the morphism of $q$-characters for $\U_q(\Lo\Glie_J)\subset\U_q(\Lo\Glie)$. 
\\For a monomial $m=\underset{i\in I, a\in\CC^*}{\prod}Y_{i,a}^{u_{i,a}(m)}$ denote $\omega(m)=\underset{i\in I,a\in\CC^*}{\sum}u_{i,a}(m)\Lambda_i$. $m$ is said to be 
$J$-dominant if for all $j\in J, a\in\CC^*$ we have $u_{j,a}(m)\geq 0$. An $I$-dominant monomials is said to be dominant. $B_J$ is the set of $J$-dominant monomials, $B$ is the set of dominant monomials. 

\noindent Note that $\chi_q, \chi_q^J$ can also be defined for finite dimensional $\U_q(\Lo \Hlie)$-modules in the same way.

\noindent In the following for $M$ a finite dimensional $\U_q(\Lo\Glie)$-module, we denote by $\mathcal{M}(M)$ the set 
of monomials of $\chi_q(M)$. 

\noindent For $i\in I, a\in\CC^*$ we set : 
\begin{equation}\label{aia}A_{i,a}=Y_{i,aq_i^{-1}}Y_{i,aq_i}\underset{j/C_{j,i}=-1}{\prod}Y_{j,a}^{-1}
\underset{j/C_{j,i}=-2}{\prod}Y_{j,aq^{-1}}^{-1}Y_{j,aq}^{-1}
\underset{j/C_{j,i}=-3}{\prod}Y_{j,aq^2}^{-1}Y_{j,a}^{-1}Y_{j,aq^{-2}}^{-1}.\end{equation}

\noindent As the $A_{i,a}^{-1}$ are algebraically independent \cite{Fre} (because $C(z)$ is invertible), for $M$ a 
product of $A_{i,a}^{-1}$ we can define $v_{i,a}(M)\geq 0$ by $M=\underset{i\in I, 
a\in\CC^*}{\prod}A_{i,a}^{-v_{i,a}(m)}$. We put $v(M)=\underset{i\in I, a\in\CC^*}{\sum}v_{i,a}(m)$. 

\noindent For $m,m'\in A$, we denote $m\leq m'$ if $m'm^{-1}$ is a product of $A_{i,a}$ ($i\in I, a\in\CC^*$).

\noindent For $\lambda\in - Q^+$ we set $v(\lambda)=-\lambda(\Lambda_1^{\vee}+...+\Lambda_n^{\vee})$. For $M$ a 
product of $A_{i,a}^{-1}$, we have $v(M)=v(\omega(\lambda))$.

\begin{defi}\label{monomrn}\cite{Fre2} A monomial $m\in A-\{1\}$ is said to be right-negative if for all $a\in\CC^*$, 
for $L=\text{max}\{l\in\ZZ/\exists i\in I, u_{i,aq^L}(m)\neq 0\}$, we have $\forall j\in I$, 
$u_{j,aq^L}(m)\neq 0\Rightarrow 
u_{j,aq^L}(m)<0$.\end{defi}

\noindent Note that a right-negative monomial is not dominant.

\begin{lem}\cite{Fre2}\label{rn} 1) For $i\in I, a\in\CC^*$, $A_{i,a}^{-1}$ is right-negative.

2) A product of right-negative monomials is right-negative.

3) If $m$ is right-negative, then $m'\leq m$ implies that $m'$ is right-negative.\end{lem}

\noindent Let $\beta : \ZZ[Y_{j,b}^{\pm}]_{j\in I,b\in\CC^*}\rightarrow \ZZ[e^{\omega}]_{\omega\in P}$ be the ring morphism such that $\beta(m)=e^{\omega(m)}$. 

\begin{prop}\label{restrcar}\cite{Fre} For a module $V\in\text{Rep}(\U_q(\Lo\Glie))$, let $\text{Res}(V)$ be the restricted 
$\U_q(\Glie)$-module. We have $(\beta\circ\chi_q)(V)=\chi(\text{Res}(V))$.\end{prop}

\subsubsection{$l$-highest weight representations}

\noindent The irreducible finite dimensional $\U_q(\Lo\Glie)$-modules have been classified by Chari-Pressley. They are parameterized by dominant monomials : 

\begin{defi}\label{lhigh} A $\U_q(\Lo\Glie)$-module $V$ is said to be of $l$-highest weight $m\in A$ if there is $v\in V_m$ such 
that $V=\U_q(\Lo\Glie)^-.v$ and $\forall i\in I, m\in\ZZ, x_{i,m}^+.v=0$.\end{defi}

\noindent For $m\in A$, there is a unique simple module $L(m)$ of $l$-highest weight $m$. 

\begin{thm}\cite{Cha2} The dimension of $L(m)$ is finite if and only if $m\in B$.\end{thm}

\begin{defi} For $i\in I, a\in\CC^*, k\geq 1$, the 
Kirillov-Reshetikhin module $W_{k,a}^{(i)}$ is the simple $\U_q(\Lo\Glie)$-modules of $l$-highest weight 
$m_{k,a}^{(i)}=\underset{s=1...k}{\prod}Y_{i,aq_i^{2s-2}}$.\end{defi}

\noindent We denote by $W_{0,a}^{(i)}$ the trivial representation (it is of dimension $1$). For $i\in I$ and $a\in\CC^*$, $W_{1,a}^{(i)}$ is called a fundamental representation and is denoted by $V_{i,a}$. 

\noindent The monomials $m_1=m_{k_1,a_1}^{(i)}$, $m_2=m_{k_2,a_2}^{(i)}$ are said to be in special position if the monomial $m_3=\underset{a\in\CC^*}{\prod}Y_{i,a}^{\text{max}(u_{i,a}(m_1),u_{i,a}(m_2))}$ is of the form $m_3=m_{k_3,a_3}^{(i)}$ and $m_3\neq m_1, m_3\neq m_2$.

\noindent A normal writing of an dominant monomial $m$ is a product decomposition $m=\underset{i=1,...,L}{\prod}m_{k_l,a_l}^{(i_l)}$ such that for $l\neq l'$, if $i_l=i_{l'}$ then $m_{k_l,a_l}^{(i_l)}$, $m_{k_{l'},a_{l'}}^{(i_{l'})}$ are not in special position. Any dominant monomial has a unique normal writing up to permuting the monomials (see \cite{Cha2}).

\noindent It follows from the study of the representations of $\U_q(\Lo sl_2)$ in \cite{Cha0, Cha, Fre} that :

\begin{prop}\label{aidesldeux} Suppose that $\Glie=sl_2$.

(1) $W_{k,a}$ is of dimension $k+1$ and :
$$\chi_q(W_{k,a})=m_{k,a}(1+A_{aq^{2k-1}}^{-1}(1+A_{aq^{2(k-1)-1}}^{-1}(1+...(1+A_{aq^{2-1}}^{-1}))...).$$

(2) $L(Y_a)\otimes L(Y_{aq^2})\otimes ...\otimes L(Y_{aq^{2(k-1)}})$ is of dimension $2^k$ and of $q$-character :
$$m_{k,a}(1+A_{aq}^{-1})(1+A_{aq^3}^{-1})...(1+A_{aq^{2k-1}}^{-1}).$$
In particular all $l$-weight spaces are of dimension $1$.

(3) for $m$ a dominant monomial and $m=m_{k_1,a_1}...m_{k_l,a_l}$ a normal writing we have :
$$L(m)\simeq W_{k_1,a_1}\otimes ...\otimes W_{k_l,a_l}.$$
\end{prop}

\subsubsection{Complementary reminders}\label{compred}

\begin{defi}\cite{Nab} A $\U_q(\Lo\Glie)$-module is said to be special if his $q$-character has a unique dominant 
monomial.\end{defi}

\noindent Note that a special module is a simple $l$-highest weight module. But in general all simple $l$-highest
weight module are not special.

\noindent For example for $\Glie=sl_2$ the Kirillov-Reshetikhin modules are special (see proposition
\ref{aidesldeux}), and :

\begin{thm}\cite{Fre2} The fundamental representations are special.\end{thm}

\begin{thm}\cite{Nab, Nad} For simply-laced cases, the Kirillov-Reshetikhin modules are special.\end{thm}

\noindent We will prove (theorem \ref{domkr}) that the Kirillov-Reshetikhin modules are special for all types.

\noindent In \cite{Fre2} an algorithm is proposed to compute $q$-characters of special $l$-highest weight modules (the proposition \ref{jdecomp} with $|J|=1$ is a formal reinterpretation of it). As a consequence :

\begin{cor}\label{lowfund}\cite{Fre2} For $m\in \mathcal{M}(V_{i,a})$, $m\neq Y_{i,a}\Rightarrow m\leq 
Y_{i,a}A_{i,aq_i}^{-1}$ and $m\in\ZZ[Y_{j,aq^l}^{\pm}]_{j\in I, l > 0}$.\end{cor}

\noindent In particular for $V$ a $l$-highest weight module of highest weight monomial $m$, for all $m'\in\mathcal{M}(m)$, we have $m'\leq m$ and the $v_{i,a}(m'm^{-1}), v(m'm^{-1})=v(\omega(m')-\omega(m))\geq 0$ are well-defined.

\noindent For $J\subset I$ and $m\in A$ denote $m^{(J)}=\underset{j\in J, a\in\CC^*}{\prod}Y_{j,a}^{u_{j,a}(m)}$. For $j\in J, a\in\CC^*$ consider $A_{j,a}^{J, \pm}=(A_{j,a}^{\pm})^{(J)}$. Define $\mu_J^I:\ZZ[A_{j,a}^{J, \pm}]_{j\in J,a\in\CC^*}\rightarrow \ZZ[A_{j,a}^{\pm}]_{j\in J,a\in\CC^*}$ the ring morphism such that $\mu_J^I(A_{j,a}^{J, \pm})=A_{j,a}^{\pm}$. For $m\in B_J$, denote $L^J(m^{(J)})$ defined for $\Glie_J$. Define :
$$L_J(m)=m^{(I-J)}\mu_J^I((m^{(J)})^{-1}L_J(m^{(J)}))$$
We have :

\begin{prop}\label{jdecomp}(\cite{her05}, proposition 3.9) For a module $V\in\text{Rep}(\U_q(\Lo \Glie))$ and $J\subset I$, there is unique
decomposition in a finite sum : 
$$\chi_q(V)=\underset{m'\in B_J}{\sum}\lambda_J(m')L_J(m').$$
Moreover for all $m'$, $\lambda_J(m')\geq 0$.\end{prop}

\noindent (In \cite{her04} the $\lambda_J(m')\geq 0$ were assumed, but the proof of the uniqueness does not depend on it).

\section{The Kirillov-Reshetikhin conjecture and $T$-system}\label{deux}

Let $i\in I, a\in\CC^*, k\geq 1$. We investigate the $q$-character of the Kirillov-Reshetikhin module $W_{k,a}^{(i)}$ 
: in this paper we prove the Kirillov-Reshetikhin conjecture and that $\chi_q(W_{k,a}^{(i)})$ satisfies the $T$-system.

\subsection{Statement of the Kirillov-Reshetikhin conjecture}

For $i\in I, k\geq 1$ consider the Kirillov-Reshetikhin module restricted to $\U_q(\Glie)$ : $Q_k^{(i)}=Res(W_{k,a}^{(i)})$ (it is independent of $a\in\CC^*$). Let $\mathcal{Q}_k^{(i)}=e^{-k\Lambda_i}\chi(Q_k^{(i)})$ be his the normalized 
character. The Kirillov-Reshetikhin conjecture is the statement of the theorems \ref{conjun} and \ref{conjdeux} proved in this paper :

\begin{defi}\label{explicit} For a sequence $\nu=(\nu_k^{(i)})_{i\in I, k> 0}$ such that for all but finitely many 
$\nu_k^{(i)}$ are non zero let us define :
$$\mathcal{F}(\nu)=\underset{N=(N_k^{(i)})}{\sum}\underset{i\in I, k>0
}{\prod}\begin{pmatrix}P_k^{(i)}(\nu,N)+N_k^{(i)}\\N_k^{(i)}\end{pmatrix}e^{-k N_k^{(i)}\alpha_i}$$ 
where 
$$P_k^{(i)}(\nu,N)=\underset{l=1...\infty}{\sum}\nu_l^{(i)}\text{min}(k,l)-\underset{j\in
I,l>0}{\sum}N_l^{(j)}r_iC_{i,j}\text{min}(k/r_j,l/r_i),$$
$$\begin{pmatrix}a\\b\end{pmatrix}=\frac{\Gamma(a+1)}{\Gamma(a-b+1)\Gamma(b+1)}.$$\end{defi}

\noindent The above formulas are called non-deformed fermionic formulas (we use here the version of \cite{ki1, ki2, hkoty, knt}; the version of \cite{kr} is  slightly different because the definition of binomial coefficients is a little changed, see \cite{knt}).

\begin{thm}[The Kirillov-Reshetikhin conjecture]\label{conjun} For a sequence $\nu=(\nu_k^{(i)})_{i\in I, k> 0}$ such that for all but finitely many 
$\nu_k^{(i)}$ are zero, we set $\mathcal{Q}_{\nu}=\underset{i\in I, k\geq 1}{\prod}(\mathcal{Q}_k^{(i)})^{\nu_k^{(i)}}$. Then we have :
$$\mathcal{Q}_{\nu}\underset{\alpha\in\Delta_+}{\prod}(1-e^{-\alpha})=\mathcal{F}(\nu).$$\end{thm}

\noindent In particular the formula obtained for $Q_{\nu}=\underset{i\in I, k\geq 1}{\prod}(\chi(Q_k^{(i)}))^{\nu_k^{(i)}}$ is the character of a $\U_q(\Glie)$-module (and so invariant by the Weyl group action). Although this formula can be given explicitly from the definition \ref{explicit}, this combinatorial consequence was conjectural for exceptional types (see \cite{k}).

\subsection{$Q$-system}

\noindent For $k\geq 1$ and $i\in I$ consider the $\U_q(\Glie)$-modules $R_k^{(i)}$ defined by :

for $r_i\geq 2$ : 
$$R_k^{(i)}=(\underset{j/C_{j,i}=-1}{\bigotimes}Q_k^{(j)})\otimes(\underset{j/C_{j,i}\leq
-2}{\bigotimes}Q_{r_ik}^{(j)}),$$

for $r_i=1$ and $\Glie$ not of type $G_2$ : 
$$k=2r\text{ : }R_k^{(i)}=(\underset{j/C_{i,j}=-1}{\bigotimes}Q_k^{(j)})\otimes 
(\underset{j/C_{i,j}=-2}{\bigotimes}Q_r^{(j)}\otimes 
Q_r^{(j)}),$$
$$k=2r+1\text{ : }R_k^{(i)}=(\underset{j/C_{i,j}=-1}{\bigotimes}Q_k^{(j)})\otimes 
(\underset{j/C_{i,j}=-2}{\bigotimes}Q_{r+1}^{(j)}\otimes  
Q_r^{(j)}),$$

for $r_i=1$ and $\Glie$ of type $G_2$ let $j\in I$ such that $j\neq i$ :
$$k=3r\text{ : }R_k^{(i)}=Q_r^{(j)}\otimes Q_r^{(j)}\otimes Q_r^{(j)},$$
$$k=3r+1\text{ : }R_k^{(i)}=Q_{r+1}^{(j)}\otimes Q_r^{(j)}\otimes Q_r^{(j)},$$
$$k=3r+2\text{ : }R_k^{(i)}=Q_{r+1}^{(j)}\otimes Q_{r+1}^{(j)}\otimes Q_r^{(j)}.$$

\begin{thm}[The $Q$-system]\label{conjdeux} Let $a\in\CC^*, k\geq 1, i\in I$. 

\noindent (1) We have :
$$Q_k^{(i)}\otimes Q_k^{(i)}=Q_{k+1}^{(i)}\otimes Q_{k-1}^{(i)}\oplus R_k^{(i)}.$$
(2) $\mathcal{Q}_k^{(i)}$ considered as a polynomial in $e^{-\alpha_j}$ has a limit as a formal 
power series :
$$\exists \underset{k\rightarrow \infty}{\text{lim}}\mathcal{Q}_k^{(i)}\in\ZZ[[e^{-\alpha_j}]]_{j\in I}.$$
\end{thm}

\noindent Note that in the simply-laced case the $Q$-system is :$Q_k^{(i)}\otimes Q_k^{(i)}=Q_{k+1}^{(i)}\otimes Q_{k-1}^{(i)}\oplus
\underset{j/C_{i,j}=-1}{\bigoplus}Q_k^{(j)}$.

\noindent (1) is the $Q$-system defined in \cite{kr, ki3}. It is proved in \cite{hkoty} that if a solution of the $Q$-system is a sum of character of $\Glie$-modules and satisfies a certain asymptotic property, then it equals the formulas $\mathcal{F}(\nu)$. In \cite{k} this asymptotic property was simplified to the property (2) (in \cite{k} it is not assumed that the solution of the $Q$-system is a sum of character of $\Glie$-modules; see the section \ref{partdeuxcomp} for more comments on the asymptotic property of \cite{hkoty}). So the theorem \ref{conjdeux} implies the theorem \ref{conjun} .

\noindent We will prove a stronger version of the theorem \ref{conjdeux} called $T$-system (theorem \ref{tsyst}). 

\subsection{$T$-system} The $T$-system was introduced in \cite{kns} as a system of functional relations associated with solvable lattice models. Motivated by results of \cite{Fre}, it was conjectured in \cite{kosy} that the
$q$-characters of Kirillov-Reshetikhin modules solve the $T$-system. This conjecture is proved in this paper
(theorem \ref{tsyst}). For simply-laced cases, the conjecture was proved by Nakajima \cite{Nab, Nad} with the 
help of $q,t$-characters, and in particular with the main result of \cite{Nab} whose 
proof involves quiver varieties. Although $q,t$-characters can be defined in general (see \cite{her02}), this result 
in \cite{Nab} has not been proved for non simply-laced cases (see the conjectures of \cite{her02}). The method of the proof used in this paper is different and based on the result of the theorem  \ref{domkr} (an induction argument of \cite{her05} (lemma \ref{submod}) is used). As it is
purely algebraic, it can be uniformly extended to non simply-laced cases.

For $i\in , k\geq 1, a\in\CC^*$ define the $\U_q(\Lo\Glie)$-module $S_{k,a}^{(i)}$ by : 

for $r_i\geq 2$ :
$$S_{k,a}^{(i)}=(\underset{j/C_{j,i}=-1}{\bigotimes}W_{k,aq_i}^{(j)})\otimes(\underset{j/C_{j,i}\leq
-2}{\bigotimes}W_{r_ik,aq}^{(j)}),$$

for $r_i=1$ and $\Glie$ not of type $G_2$ :
$$k=2r\text{ : }S_{k,a}^{(i)}=(\underset{j/C_{i,j}=-1}{\bigotimes}W_{k,aq}^{(j)})\otimes
(\underset{j/C_{i,j}=-2}{\bigotimes}W_{r,aq}^{(j)}\otimes
W_{r,aq^3}^{(j)}),$$
$$k=2r+1\text{ : }S_{k,a}^{(i)}=(\underset{j/C_{i,j}=-1}{\bigotimes}W_{k,aq}^{(j)})\otimes
(\underset{j/C_{i,j}=-2}{\bigotimes}W_{r+1,aq}^{(j)}\otimes W_{r,aq^3}^{(j)}),$$

for $r_i=1$ and $\Glie$ of type $G_2$ let $j\in I$ such that $j\neq i$ :
$$k=3r\text{ : }S_{k,a}^{(i)}=W_{r,aq}^{(j)}\otimes W_{r,aq^3}^{(j)}\otimes W_{r,aq^5}^{(j)},$$
$$k=3r+1\text{ : }S_{k,a}^{(i)}=W_{r+1,aq}^{(j)}\otimes W_{r,aq^3}^{(j)}\otimes W_{r,aq^5}^{(j)},$$
$$k=3r+2\text{ : }S_{k,a}^{(i)}=W_{r+1,aq}^{(j)}\otimes W_{r+1,aq^3}^{(j)}\otimes W_{r,aq^5}^{(j)}.$$

\noindent Remark : we will see later, in view of the lemma \ref{commute} and the proposition \ref{xspecial}, that in all cases the tensor products of the modules involved in the definition of $S_{k,a}^{(i)}$ commute, and so $S_{k,a}^{(i)}$ is well-defined. However until section \ref{sequ}, we
only consider $\chi_q(S_{k,a}^{(i)})$ which is clearly well-defined.

\begin{thm}[The $T$-system]\label{tsyst} Let $a\in\CC^*, k\geq 1, i\in I$. 

\noindent (1) We have :
$$\chi_q(W_{k,a}^{(i)})\chi_q(W_{k,aq_i^2}^{(i)})=\chi_q(W_{k+1,a}^{(i)})\chi_q(W_{k-1,aq_i^2}^{(i)})
+\chi_q(S_{k,a}^{(i)}).$$ 
(2) The normalized $q$-character of $W_{k,a}^{(i)}$ considered as a polynomial in $A_{j,b}^{-1}$ has a limit as a
formal power series : 
$$\exists \underset{k\rightarrow
\infty}{\text{lim}}\frac{\chi_q(W_{k,aq_i^{-2k}}^{(i)})}{m_{k,aq_i^{-2k}}^{(i)}}\in\ZZ[[A_{j,aq^m}^{-1}]]_{j\in I, m\in\ZZ}.$$\end{thm}

\noindent Note that in the simply-laced cases the T-system is :
$$\chi_q(W_{k,a}^{(i)})\chi_q(W_{k,aq^2}^{(i)})=\chi_q(W_{k+1,a}^{(i)})\chi_q(W_{k-1,aq^2}^{(i)})
+\underset{j/C_{i,j}=-1}{\prod}\chi_q(W_{k,aq}^{(i)})$$
(see the section \ref{exptsyst} for the other formulas of the $T$-systems)

\noindent The theorem \ref{tsyst} implies the theorem \ref{conjdeux} because $\text{Res}(W_{k,a}^{(i)})=Q_k^{(i)}$ and
$\text{Res}(S_{k,a}^{(i)})=R_k^{(i)}$. The theorem \ref{tsyst} is proved in section \ref{proofth} with the main result of the section \ref{special}. 

\noindent Note that the methods of our proof may also be used to prove the conjectures of \cite{hkott} for twisted cases, and to
establish a $T$-system for representations of general quantum affinizations (such as quantum toroidal algebras) considered in \cite{her04}. We let these points for another paper.

\noindent First we need the following result :

\section{The Kirillov-Reshetikhin modules are special}\label{special} In this section we prove :

\begin{thm}\label{domkr} The Kirillov-Reshetikhin modules are special.\end{thm}

\noindent This result was proved for fundamental representations in \cite{Fre2}, for Kirillov-Reshetikhin modules of 
type $sl_2$ in \cite{Fre} (see proposition \ref{aidesldeux}), and in simply-laced cases in \cite{Nab, Nad}. An alternative 
proof for fundamental representations was proposed in \cite{her05}. Arguments of this last proof are used here.

\noindent This result implies that the Frenkel-Mukhin algorithm works for Kirillov-Reshetikhin modules. In particular we could have informations on the structure of their $q$-character (for example we will prove in a paper in preparation precise results for $l$-weight spaces of Kirillov-Reshetikhin modules of type $A,B$ with technics developed in \cite{her05}).

\subsection{Preliminary results} The following ingredient is used in the proof :

\begin{lem}\label{submod}(\cite{her05}, lemma 3.3, 3.4) Let $V$ be a finite dimensional $\U_q(\Lo\Glie)$-module .

\noindent (i) For $W\subset V$ a $\U_q(\Lo\Hlie)$-submodule of $V$ and $i\in I$,
$W_i'=\underset{r\in\ZZ}{\sum}x_{i,r}^-.W$ is a $\U_q(\Lo\Hlie)$-submodule of $V$.

\noindent (ii) Suppose that $\Glie=sl_2$. For $p\in\ZZ$ let $L_{\geq p}=\underset{q \geq p}{\sum}L_{q\Lambda}$ and 
$L'_{\geq p}=\underset{r\in\ZZ}{\sum}x_r^-.L_{\geq p}$. Then $L_{\geq p}, L_{\geq p}'$ are $\U_q(\Lo\Hlie)$-submodule of $L$ and
$(L_{\geq p}')_{m}\neq 0\Rightarrow \exists m'$, $(L_{\geq p})_{m'}\neq \{0\}$ and $m \leq m'$.\end{lem}

\noindent We will also use a result of \cite{Fre2} in the more precise form of \cite{her04} :

\noindent Let $i\in I$, $\Hlie_i^{\perp}=\{\mu\in\Hlie/\alpha_i(\mu)=0\}$ and $A^{(i)}$ be the commutative group 
of monomials generated by variables $Y_{i,a}^{\pm}$ ($a\in\CC^*$), $k_{\mu}$ ($\mu\in\Hlie_i^{\perp}$), 
$Z_{j,c}^{\pm}$ ($j\neq i$, $c\in\CC^*$). Let $\tau_i:A\rightarrow A^{(i)}$ be the group morphism defined by ($j\in I$, 
$a\in\CC^*$): 
$$\tau_i(Y_{j,a})=Y_{j,a}^{\delta_{j,i}}\underset{k\neq 
i}{\prod}\underset{r\in\ZZ}{\prod}Z_{k,a q^r}^{p_{j,k}(r)}k_{\nu(\Lambda_j)-\delta_{j,i}r_i\alpha_i^{\vee}/2}.$$ 
The $p_{i,j}(r)\in\ZZ$ are defined in the following way : we write 
$\tilde{C}(z)=\frac{\tilde{C}'(z)}{d(z)}$ where $d(z), \tilde{C}'_{i,j}(z)\in\ZZ[z^{\pm}]$ and 
$(D(z)\tilde{C}'(z))_{i,j}=\underset{r\in\ZZ}{\sum}p_{i,j}(r)z^r$. Note that we have $\nu(\Lambda_j)-\delta_{j,i}r_i\alpha_i^{\vee}/2\in \Hlie_i^{\perp}$ because $\alpha_i(\nu(\Lambda_j)-\delta_{j,i}r_i\alpha_i^{\vee}/2)=\Lambda_j(r_i\alpha_i^{\vee})-r_i\delta_{i,j}=0$.

\noindent This morphism was defined in \cite{Fre2}, section 3.3 without the terms $k_{\mu}$, and was then refined in \cite{her04}, section 5.5.2 with this term (which will be used in the section \ref{pf}).

\noindent For $M\in A^{(i)}$ consider $\mu(M)\in\Hlie_i^{\perp}$, $u_{i,a}(M)\in\ZZ$, such that : 
$$M\in k_{\mu(M)}\underset{a\in\CC^*}{\prod}Y_{i,a}^{u_{i,a}(M)}\ZZ[Z_{j,c}^{\pm}]_{j\neq i, c\in\CC^*}.$$
We also set $u_i(M)=\underset{a\in\CC^*}{\sum}u_{i,a}(M)$. Note that for $m\in A$ and $a\in\CC^*$ we have $u_{i,a}(m)=u_{i,a}(\tau_i(m))$ and :
$$\nu(\omega(m))=\mu(\tau_i(m))+u_i(m)r_i\alpha_i^{\vee}/2=\mu(\tau_i(m))+u_i(\tau_i(m))r_i\alpha_i^{\vee}/2,$$
or equivalently :
$$\mu(\tau_i(m))=\nu(\omega(m))-\alpha_i(\nu(\omega(m)))\alpha_i^{\vee}/2$$
(see the definition of \cite{her04}, section 5.5.2).

\begin{lem}\label{aidedeux} Let $V\in\text{Rep}(\U_q(\Lo\Glie))$ and consider a decomposition $\tau_i(\chi_q(V))=\underset{r}{\sum}P_rQ_r$ where $P_r\in\ZZ[Y_{i,a}^{\pm}]_{a\in\CC^*}$,
$Q_r$ is a monomial in $\ZZ[Z_{j,c}^{\pm}, k_{\lambda}]_{j\neq i,c\in\CC^*, \lambda\in\Hlie_i^{\perp}}$ and all 
monomials $Q_r$ are distinct. Then the
$\U_q(\Lo\Glie_i)$-module $V$ is isomorphic to a direct sum $\underset{r}{\bigoplus}V_r$ where
$\chi_q^i(V_r)=P_r$.\end{lem}

\noindent This result was proved in \cite{Fre2}, lemma 3.4 without the term $k_{\mu}$, and in \cite{her04}, lemma 5.10 the proof was extended for the terms $k_{\mu}$.

\subsection{Proof of the theorem \ref{domkr}}\label{pf} The theorem \ref{domkr} is a consequence of :

\begin{lem}\label{domkrlem} For $m\in\mathcal{M}(W_{k,a}^{(i)})$, we have $m\neq m_{k,a}^{(i)}\Rightarrow m\leq 
m_{k,a}^{(i)}A_{i,aq_i^{2k-1}}^{-1}$. In particular $m$ is right-negative and not dominant.\end{lem}

\demo For $m\leq m_{k,a}^{(i)}$ we denote $w(m)=v(m(m_{k,a}^{(i)})^{-1})$. Let $M=V_{i,a}\otimes V_{i,aq_i^2}\otimes ...\otimes 
V_{i,aq_i^{2k-2}}$. $m_{k,a}^{(i)}$ is the monomial of highest weight in 
$\mathcal{M}(M)$. $M$ is an $l$-highest weight module (see \cite{c, kas}; 
we do 
not really need this point as we could work with $M'\subset M$ the 
submodule generated by an $l$-highest weight vector $v$). So $W_{k,a}^{(i)}$ is 
the quotient of $M$ by the maximal sub $\U_q(\Lo\Glie)$-module $N\subset 
M$.

\noindent Consider the sub-$\U_q(\Lo\Glie_i)$-submodule $M_i'$ of $M$ generated by an highest weight vector $v$. It is an $l$-highest weight
$\U_q(\Lo\Glie_i)$-module of $l$-highest weight $m_{k,a}^{(i)}$, and so it has a simple quotient which is a
Kirillov-Reshetikhin module $L_i$ of type $sl_2$. So there is $N_i$ a maximal $\U_q(\Lo\Glie_i)$-submodule of $M_i'$ such that $(M_i'/N_i)\simeq L_i$. For $j\neq i$ and $R\geq 0$, we have $M_{k\Lambda_i-R\alpha_i+\alpha_j}=\{0\}$ and so for all $m\in\ZZ$, we have $x_{j,m}^+.M_i'=\{0\}$. So $\U_q(\Lo\Glie).N_i$ is a proper submodule of
$M'$. So $W_{k,a}^{(i)}$ is subquotient of $M/\U_q(\Lo\Glie).N_i$. In particular the
$\U_q(\Lo\Glie_i)$-submodule $M_i$ of $W_{k,a}^{(i)}$ generated by $v$ is simple and isomorphic to $L_i$.

\noindent It follows from the corollary \ref{lowfund} that :
$$\underset{m\leq m_{k,a}^{(i)}/w(m)=1}{\sum}M_m\subset M_{k\Lambda_i^{\vee}-\alpha_i}\subset
\U_q(\Lo\Glie_i).v.$$
In particular $\underset{m\leq m_{k,a}^{(i)}/w(m)=1}{\sum}(W_{k,a}^{(i)})_m\subset M_i$. So it follows from the last
paragraph and from the (1) of the proposition \ref{aidesldeux} that :
$$\underset{m\leq m_{k,a}^{(i)}/w(m)=1}{\sum}(W_{k,a}^{(i)})_m=(W_{k,a}^{(i)})_{m_{k,a}^{(i)}A_{i,aq_i^{2k-1}}^{-1}}$$ 
and that this space is of dimension $1$.

\noindent Now consider $m\in\mathcal{M}(W_{k,a}^{(i)})$ such that $m\neq m_{k,a}^{(i)}$, and let us prove by induction on $w(m)\geq 1$ that $m\leq 
m_{k,a}^{(i)}A_{i,aq_i^{2k-1}}^{-1}$. For $w(m)=1$ we have proved that $m=m_{k,a}^{(i)}A_{i,aq_i^{2k-1}}^{-1}$. In general suppose
that $w(m)=p+1$ ($p\geq 1$). It follows from the structure of $M_i\simeq L_i$ (which is also a $\U_q(\Lo\Hlie)$-module) that we can suppose that $(M_i)_m=\{0\}$. Consider :
$$W=\underset{m'\leq m_{k,a}^{(i)}/w(m')\leq 
p}{\bigoplus}(W_{k,a}^{(i)})_{m'}=\underset{\lambda\in\Hlie/v(\lambda-k\Lambda_i)\leq p}{\bigoplus}(W_{k,a}^{(i)})_{\lambda}.$$ 
Note that $W$ is a $\U_q(\Lo\Hlie)$-submodule of $W_{k,a}^{(i)}$. As $W_{k,a}^{(i)}$ is a $l$-highest weight modules, we have :
$$\underset{m'\leq m_{k,a}^{(i)}/w(m')=p+1}{\bigoplus}(W_{k,a}^{(i)})_{m'}=\underset{\lambda\in\Hlie/v(\lambda-k\Lambda_i)=p+1}{\bigoplus}(W_{k,a}^{(i)})_{\lambda}\subset \underset{j\in I}{\sum}W_j\text{ where }W_j=\underset{r\in\ZZ}{\sum}x_{j,r}^-.W.$$ 
For $j\in I$, $W_j$ is
a $\U_q(\Lo\Hlie)$-submodule of $W_{k,a}^{(i)}$ ((i) of lemma \ref{submod}). So $\exists j\in I$, $(W_j)_m\neq 
\{0\}$.

\noindent Consider the decomposition $\tau_j(\chi_q(W_{k,a}^{(i)}))=\underset{r}{\sum}P_rQ_r$ of the lemma 
\ref{aidedeux} and the decomposition of $W_{k,a}^{(i)}$ as a $\U_q(\Lo\Glie_j)$-module:  
$W_{k,a}^{(i)}=\underset{r}{\bigoplus}V_r$. For a given $r$, consider $M_r\in\mathcal{M}(W_{k,a}^{(i)})$ such that 
$\tau_j(M_r)$ appears in $P_rQ_r$. For another such $M$, we have $\mu(\tau_j(M))=\mu(\tau_j(M_r))$ and so 
$\omega(MM_r^{-1})=u_j(\tau_j(MM_r^{-1}))\alpha_j^{\vee}/2$, and : 
$$u_j(\tau_j(M))=u_j(\tau_j(M_r))-2w(M)+2w(M_r)=2(p-w(M))+p_r$$ where $p_r=-2p+2w(M_r)+u_j(\tau_j(M_r))$ (it does not 
depend of $M$). So we have $w(M)\leq p\Leftrightarrow u_j(\tau_j(M))\geq p_r$. So 
$W=\underset{r}{\bigoplus}((V_r)_{\geq p_r})=\underset{r}{\bigoplus}(V_r\cap W)$. As the $V_r$ are sub 
$\U_q(\Lo\Glie_j)$-modules of $W_{k,a}^{(i)}$, we have $W_j=\underset{r}{\bigoplus}(V_r\cap W_j)$. \\Let $R$ such that 
$\tau_j(m)$ is a monomial of $P_RQ_R$. We can apply the (ii) of lemma \ref{submod} to the $\U_{q}(\Lo\Glie_j)$-module 
$V_R$ with $p_R$ for $Q_R^{-1}\tau_j(m)$ : we get that there is $M'$ a monomial of $\chi_q^j(W)$ such that 
$Q_R^{-1}\tau_i(m)\in M' \ZZ[(Y_{j,a}Y_{j,aq_j^2})^{-1}]_{a\in\CC^*}$. Consider $m'=\tau_j^{-1}(Q_RM')$ (it is a 
monomial of $\chi_q(W)$). Note that we can suppose that $m'\neq m_{k,a}^{(i)}$ (if $m'= m_{k,a}^{(i)}$, consider 
$Z=\U_q(\Lo\Glie_j).(W_{k,a}^{(i)})_{m_{k,a}^{(i)}}=(W_{k,a}^{(i)})_{m_{k,a}^{(i)}}$ or $M_i$. As $Z$ is a sub 
$\U_q(\Lo\Glie_j)$-module of $V_R$ and $Q_R^{-1}\tau_j(m)$ is not a monomial of $\chi_q^j(Z)$, we can use above 
$V_R/Z$ instead of $V_R$).

\noindent It follows from \cite{Fre2}, lemma 3.5 that 
$\tau_j(A_{j,aq_j})=Y_{j,a}Y_{j,aq_j^2}k_0$ (see \cite{her04}, lemma 5.9 for the term $k_0$). So $m < m'$ and $m\in m'\ZZ[A_{j,b}^{-1}]_{b\in\CC^*}$. As $m'\neq m_{k,a}^{(i)}$ we have $w(m')\geq 1$ and the induction hypothesis gives the result.\qed

\section{Proof of the theorem \ref{tsyst}}\label{proofth}

\subsection{Preliminary results}\label{prel} First let us prove the following :

\begin{lem}\label{commute} Let $V$ be a special module. Suppose that $V\simeq V_1\otimes
....\otimes V_r$ where $V_1,...,V_r$ are $l$-highest weight modules. Then $V_1,...,V_r$ are special and for all
$\sigma$ permutation of $\{1,...,r\}$ we have $V\simeq V_{\sigma(1)}\otimes ...\otimes V_{\sigma(r)}$.\end{lem}

\demo Let $m_1,...,m_r$ be the monomials of highest weight of $V_1,...,V_r$. If for $r'\in\{1,...,r\}$ a module $V_{r'}$ 
is not special, let $m_{r'}'\neq m_{r'}$ be a dominant monomial of $\mathcal{M}(V_{r'})$. Then 
$$m_1...m_{r'-1}m_{r'}'m_{r'+1}...m_r\in\mathcal{M}(V)$$ 
is dominant and not equal to the highest weight monomial $m_1...m_r$. In 
particular for all $\sigma$, $V_{\sigma(1)}\otimes ...\otimes V_{\sigma(r)}$ is special and so is simple isomorphic to 
the simple $l$-highest weight module $L(m_1...m_r)\simeq V$.\qed

\noindent Let $i\in I, k\geq 1, a\in\CC^*$. Let $M=m_{k,a}^{(i)}m_{k,aq_i^2}^{(i)}=m_{k+1,a}^{(i)}m_{k-1,aq_i^2}^{(i)}$, and 
$M'=MA_{i,aq_i^{2k-1}}^{-1}...A_{i,aq_i}^{-1}$ the highest weight monomial of $S_{k,a}^{(i)}$. Let us write the dominant monomial $M'$ in a normal way : 
\begin{equation}\label{normd}M'=\underset{l=1...L}{\prod}m_{k_l,a_l}^{(i_l)}\text{ where $i_l\neq i$}.\end{equation} 
Consider the sets of monomials : 
$$\mathcal{B}=\{m_{k,a}^{(i)}A_{i,aq_i^{2k-1}}^{-1}...A_{i,aq_i^{2(k-k')-1}}^{-1}A_{i_l,a_lq_{i_l}^{2k_l-1}}^{-1}/0\leq 
k'\leq k-1\text{ , }1\leq l\leq L\},$$
$$\mathcal{B'}=\{m_{k,a}^{(i)}, m_{k,a}^{(i)}A_{i,aq_i^{2k-1}}^{-1}, 
m_{k,a}^{(i)}A_{i,aq_i^{2k-1}}^{-1}A_{i,aq_i^{2k-3}}^{-1},...,m_{k,a}^{(i)}A_{i,aq_i^{2k-1}}^{-1}...A_{i,aq_i}^{-1}\}.$$

\begin{lem}\label{check} The monomials of $m_{k,aq_i^2}^{(i)}.\mathcal{B}$ are right negative.\end{lem}

\demo For $a\in\CC^*$ and $m\in A$, let us define $\mu_a(m)=\text{max}\{l\in\ZZ/\exists i\in I, u_{i,aq^l}(m)\neq
0\}$. Let $\alpha=MA_{i,aq_i^{2k-1}}^{-1}...A_{i,aq_i^{2(k-k')-1}}^{-1}$, and $\alpha A_{i_l,aq_{i_l}^{2k_l-1}}^{-1}\in m_{k,aq_i^2}^{(i)}\mathcal{B}$. It suffices to check that 
$\mu_a(\alpha)<\mu_a(A_{i_l,aq_{i_l}^{2k_l-1}}^{-1})$ :

Case 1: $r_i=1$ : $\mu_a(\alpha)\leq 2k-1$.

\noindent If $r_{i_l}=1$ : $a_l=aq^{2k-1}$, $\mu_a(A_{i_l,a_lq^{r_{i_l}}}^{-1})=2k+1$.

\noindent If $r_{i_l}=2$ : $a_l=aq^{2k-1}$ or $aq^{2k-3}$, $\mu_a(A_{i_l,a_lq^{r_{i_l}}}^{-1})=2k+3$ or $2k+1$.

\noindent If $r_{i_l}=3$ : $a_l=aq^{2k-1}$ or $aq^{2k-3}$ or $aq^{2k-5}$, $\mu_a(A_{i_l,a_lq^{r_{i_l}}}^{-1})=2k+5$ or 
$2k+3$ or $2k+1$.

Case 2 : $r_i=2$ : $\mu_a(\alpha)\leq 4k-1$.

\noindent If $r_{i_l}=1$ : $a_l=aq^{4k-1}$, $\mu_a(A_{i_l,a_lq^{r_{i_l}}}^{-1})=4k+1$.

\noindent If $r_{i_l}=2$ : $a_l=aq_i^{4k-1}$ or $aq_i^{4k-3}$, $\mu_a(A_{i_l,a_lq^{r_{i_l}}}^{-1})=4k+3$ or $4k+1$.

Case 3 : $r_i=3$ : $\mu_a(\alpha)\leq 6k-1$.

\noindent If $r_{i_l}=1$ : $a_l=aq^{6k-1}$, $\mu_a(A_{i_l,a_lq^{r_{i_l}}}^{-1})=6k+1$.\qed

\begin{prop}\label{xspecial} For $i\in I, k\geq 1, a\in\CC^*$, the module $S_{k,a}^{(i)}$ is special. In particular $M'$ is the unique dominant monomial of $\chi_q(S_{k,a}^{(i)})$.\end{prop}

\noindent Note that this result, combined with lemma \ref{commute}, implies that the modules in the tensor product 
$S_{k,a}^{(i)}$ commute.

\demo Let us write it : $$\chi_q(S_{k,a}^{(i)})=\chi_q(W_{k_1,a_1}^{(i_1)})... \chi_q(W_{k_L,a_L}^{(i_L)}).$$ The
monomials of $\mathcal{M}(S_{k,a}^{(i)})-\{m_{k_1,a_1}^{(i_1)}...m_{k_L,a_L}^{(i_L)}\}$ are lower than one of the
following monomials $m_{k_l,a_l}^{(i_l)}A_{i_l,a_lq_{i_l}^{2k_l-1}}^{-1}\underset{l'\neq
l}{\prod}m_{k_{l'},a_{l'}}^{(i_{l'})}$ (lemma \ref{domkrlem}). But in each case these monomials are right-negative
(these monomials are exactly the monomials of $m_{k,aq_i^2}^{(i)}\mathcal{B}$ with $k'=k-1$, and it is checked in lemma \ref{check} that these monomials are right-negative). So the monomials of
$\mathcal{M}(S_{k,a}^{(i)})-\{m_{k_1,a_1}^{(i_1)}...m_{k_L,a_L}^{(i_L)}\}$, are right-negative so not dominant.\qed

\subsection{Proof of the theorem \ref{tsyst} (1)}\label{partun}

\noindent The screening operators $S_i$ were defined in \cite{Fre}.

\begin{thm}\label{intker}\cite{Fre, Fre2} We have $\text{Im}(\chi_q)=\underset{i\in
I}{\bigcap}\text{Ker}(S_i)$. In particular a non zero element in
$\text{Im}(\chi_q)$ has at least one dominant monomial.\end{thm}

\noindent In this paper we will only use the second part of this theorem, and so we do not directly use the screening 
operators.

\noindent The two terms of the equality of the theorem \ref{tsyst} are in $\text{Im}(\chi_q)$ and so are
characterized by the coefficient of
their dominant monomials. So it suffices to determine the dominant
monomials of each product.

\noindent First let us prove the following lemma concerning the monomials of $\chi_q(W_{k,a}^{(i)})$ : 

\begin{lem}\label{moreinfo} The monomials $m$ of $\chi_q(W_{k,a}^{(i)})$ are lower than a monomial of 
$\mathcal{B}$ or are in $\mathcal{B'}$.\end{lem}

\noindent (An analog result is proved in \cite{Nad} for the simply-laced cases).

\demo We prove this statement by induction on $w(m)=v(m(m_{k,a}^{(i)})^{-1})\geq 0$. For $w(m)=0$ we have 
$m=m_{k,a}^{(i)}\in\mathcal{B'}$. For $w(m)\geq 1$ it follows from the theorem \ref{domkr} that there is 
$j\in I$ such that $m\notin B_j$. So we get from the proposition \ref{jdecomp} a monomial 
$m'\in\mathcal{M}(\chi_q(W_{k,a}^{(i)}))$ such that $w(m')<w(m)$ and $m$ is a monomial of $L_j(m')$ and $L_j(m')$ 
appears in the decomposition of $\chi_q(W_{k,a}^{(i)})$. In particular $m\leq m'$, and if $m'$ is lower than a 
monomial in $\mathcal{B}$, so is $m$. So we can suppose that $m'\in\mathcal{B'}$. If 
$m'=m_{k,a}^{(i)}$, we have $j=i$ and the monomials of $L_i(m)$ are the monomials of $\mathcal{B'}$ (see 
the proof of lemma \ref{domkrlem}, $M_i\simeq L_i$ and the proposition \ref{aidesldeux} (1)). For $m'\neq 
m_{k,a}^{(i)}$, $L_{j}(m)$ corresponds to the $q$-character of a tensor product of Kirillov-Reshetikhin modules of 
type $sl_2$ ((3) of the proposition \ref{aidesldeux}). Let $m_{k_{l'},a_{l'}}$ be the corresponding monomials (that it 
to say a normal form). For each $l'$, $a_{l'}q_{j}^{2(k_{l'}-1)}$ is equal to one $a_{l}q_{i_l}^{2(k_l-1)}$ with 
$i_l=j$ (see the decomposition (\ref{normd}) of the section \ref{prel}). We can conclude with (1) of proposition \ref{aidesldeux}.\qed

\begin{lem}\label{premier} 1) The dominant monomials of $\chi_q(W_{k,a}^{(i)})\chi_q(W_{k,aq_i^2}^{(i)})$ are 
$$M, MA_{i,aq_i^{2k-1}}^{-1},
MA_{i,aq_i^{2k-1}}^{-1}A_{i,aq_i^{2k-3}}^{-1},...,MA_{i,aq_i^{2k-1}}^{-1}...A_{i,aq_i}^{-1}.$$
2) The dominant monomials of $\chi_q(W_{k+1,a}^{(i)})\chi_q(W_{k-1,aq_i^2}^{(i)})$ are
$$M, MA_{i,aq_i^{2k-1}}^{-1}, 
MA_{i,aq_i^{2k-1}}^{-1}A_{i,aq_i^{2k-3}}^{-1},...,MA_{i,aq_i^{2k-1}}^{-1}...A_{i,aq_i^3}^{-1}.$$
In each case the dominant monomials appear with multiplicity $1$.\end{lem}

\demo We prove 1) (the proof is analog for 2)). Let $m_1\in\mathcal{M}(W_{k,a}^{(i)})$,
$m_2\in\mathcal{M}(W_{k,aq_i^2}^{(i)})$ and suppose that $m_1m_2$ is dominant. If $m_1\neq m_{k,a}^{(i)}$ and $m_2\neq
m_{k,aq_i^2}^{(i)}$, the theorem \ref{domkr} gives that $m_1$ and $m_2$ are right-negative, so $m_1m_2$ is right-negative and not dominant. If $m_2\neq m_{k,aq_i^2}^{(i)}$ we have $m_1=m_{k,a}^{(i)}$ and it follows from the lemma
\ref{domkrlem} that we have $m_2\leq m_{k,aq_i^2}^{(i)}A_{i,aq_i^3}^{-1}$, and so $m_1m_2\leq
m_{k,a}^{(i)}m_{k,aq_i^2}^{(i)}A_{i,aq_i^3}^{-1}$. But this last monomial is right-negative, so $m_1m_2$ is not
dominant. So $m_2=m_{k,aq_i^2}^{(i)}$.

\noindent Consider the decomposition (\ref{normd}) of $M'$ (section \ref{prel}). It follows from the lemma \ref{moreinfo} that the 
monomials $m$ of $\chi_q(W_{k,a}^{(i)})$ not in $\mathcal{B'}$ are lower than a monomial in 
$\mathcal{B}$. We can conclude because the monomials in $m_{k,aq_i^2}^{(i)}\mathcal{B}$ are 
right-negative (see lemma \ref{check}).\qed

{\it End of the proof of the theorem \ref{tsyst} (1) :}

\noindent The unique dominant monomial that appears in 
$\chi_q(W_{k,a}^{(i)})\chi_q(W_{k,aq_i^2}^{(i)})-\chi_q(W_{k+1,a}^{(i)})\chi_q(W_{k-1,aq_i^2}^{(i)})$ is 
$M'$, and it has a multiplicity $1$. We can conclude with the theorem \ref{intker} because $M'$ is the unique dominant monomial of $\chi_q(S_{k,a}^{(i)})$ (proposition \ref{xspecial}).\qed

\subsection{Proof of the theorem \ref{tsyst} (2)}\label{partdeux}

\subsubsection{Preliminary}

First let us see that :

\begin{lem}\label{extl} Let $a\in\CC^*$ and $m$ be a monomial in $\ZZ[Y_{j,aq^l}^{\pm}]_{j\in I, l\geq 0}$. Let $i\in 
I$ and suppose that $u_{i,a}(m)\geq 1$. Let $j\in I$ such that $m$ is $j$-dominant.

1) If  $j\neq i$ we have $Y_{i,a}L_j(Y_{i,a}^{-1}m)=L_j(m)$.

2) If $j=i$ let $n_{m'}\geq 0$ such that $L_i(m)=\underset{m'}{\sum}n_{m'}m'$. We have :
$$Y_{i,a}L_j(Y_{i,a}^{-1}m)=\underset{m'\leq m/v_{i,aq_i}(m'm^{-1})=0}{\sum}n_{m'}m'.$$\end{lem}

\demo It follows from the definition of the $L_j(M)$ (section \ref{compred}) that it suffices to look at the $sl_2$-case. But for $\Glie=sl_2$ one can use the explicit description of $L_i(m)$ in the proposition \ref{aidesldeux}.\qed

\noindent In \cite{Nad} an argument based on the Frenkel-Mukhin algorithm was used to prove the (2) of the theorem \ref{tsyst} (and the lemma \ref{aideconv}) in 
simply-laced cases. For the general case we use a different proof based on an explicit formulation of the Frenkel-Mukhin algorithm (the proposition \ref{jdecomp}). Let us prove the following lemma :

\begin{lem}\label{aideconv} We have :
$$\chi_q(W_{k+1,a}^{(i)})=Y_{i,a}\chi_q(W_{k,aq_i^2}^{(i)})+E$$
where $E\in m_{k+1,a}^{(i)}A_{i,aq_i}^{-1}A_{i,aq_i^3}^{-1}...A_{i,aq_i^{2k+1}}^{-1}\ZZ[A_{j,aq^m}^{-1}]_{j\in I, m\geq 0}$. Moreover $E$ has positive coefficients.\end{lem}

\demo First let us prove by induction on $w'(m)=v(m(m_{k,aq_i^2}^{(i)})^{-1})\geq 0$ that for 
$m\in\mathcal{M}(W_{k,aq_i^2}^{(i)})$ we have $Y_{i,a}m\in\mathcal{M}(W_{k+1,a}^{(i)})$ and the coefficient of $m$ in 
$\chi_q(W_{k,aq_i^2}^{(i)})$ is equal to the coefficient of $Y_{i,a}m$ in $\chi_q(W_{k+1,a}^{(i)})$. For 
$m=m_{k,aq_i^2}^{(i)}=Y_{i,a}^{-1}m_{k+1,a}^{(i)}$, it is clear. For $m<m_{k,aq_i^2}^{(i)}$ it follows from the theorem
\ref{domkr} that there is $j\in I$ such that $m\notin B_j$. From the proposition \ref{jdecomp} there is
$m'\in B_j\cap\mathcal{M}(W_{k,aq_i^2}^{(i)})$ such that $w'(m')<w'(m)$, $m$ is a monomial of $L_j(m')$, and $L_j(m')$
appears in the decomposition of the proposition \ref{jdecomp}. Note that the corollary \ref{lowfund} implies that we can use the lemma 
\ref{extl} for all such $m'$. It gives that the coefficient of $m$ in $L_j(m')$ is equal to the coefficient of $mY_{i,a}$ in 
$L_j(m'Y_{i,a})$. But by the induction hypothesis, the coefficients of $L_j(m')$ in $\chi_q(W_{k,aq_i^2}^{(i)})$ is equal to the 
coefficient of $L_j(m'Y_{i,a})$ in $\chi_q(W_{k+1,a}^{(i)})$ and so we get the result for $m$.

\noindent So we have proved that $E$ has positive coefficients. Consider the following property $P(m)$ of a monomial 
$m$ : $$P(m)\text{ : ``}m\in 
m_{k+1,a}^{(i)}A_{i,aq_i}^{-1}A_{i,aq_i^3}^{-1}...A_{i,aq_i^{2k+1}}^{-1}\ZZ[A_{j,aq^m}^{-1}]_{j\in I, m\geq 0}.\text{''}$$ 
To 
conclude our proof it suffices to show that for $m\in\mathcal{M}(W_{k+1,a}^{(i)})$ we have 
$m\in Y_{i,a}\mathcal{M}(W_{k,aq_i^2}^{(i)})$ or $P(m)$ is satisfied. We use an induction and we have to prove a 
little 
more in 
this induction :

\noindent For $m\in\mathcal{M}(W_{k+1,a}^{(i)})$, put $w_{j,b}(m)=v_{j,b}(m(m_{k+1,a}^{(i)})^{-1})$ and 
$w(m)=v(m(m_{k+1,a}^{(i)})^{-1})$. We prove by induction on $w(m)\geq 0$ that a monomial 
$m\in\mathcal{M}(W_{k+1,a}^{(i)})$ satisfies the property $P(m)$, or the following properties $\alpha_1(m), 
\alpha_2(m), \alpha_3(m)$ are simultaneously satisfied :

$\alpha_1(m)$ : ``$m\in Y_{i,a}\mathcal{M}(W_{k,aq_i^2}^{(i)})$.''

$\alpha_2(m)$ : ``for $k(m)=\text{max}\{k'\leq k/w_{i,aq_i^{1+2k'}}(m)=0\}$ we have :
$$w_{i,aq_i^{3+2k(m)}}(m), w_{i,aq_i^{5+2k(m)}}(m),...,w_{i,aq_i^{2k+1}}(m)\geq 1.\text{''}$$

$\alpha_3(m)$ : ``for all $j\in I$, all $l < 2r_i(k(m)+1)$ we have $w_{j,aq^{l+r_j}}(m)=0$.''

\noindent For $m=m_{k+1,a}^{(i)}=Y_{i,a}m_{k,aq_i^2}^{(i)}$, $\alpha_1(m),\alpha_3(m)$ are clear, and $\alpha_2(m)$ is 
satisfied with  
$k(m)=k$. For $m<m_{k+1,a}^{(i)}$ it follows from the theorem
\ref{domkr} that there is $j\in I$ such that $m\notin B_j$. From the proposition \ref{jdecomp} there is
$m'\in B_j\cap\mathcal{M}(W_{k+1,a}^{(i)})$ such that $w(m')<w(m)$, $m$ is a monomial of $L_j(m')$ and $L_j(m')$ 
appears in the decomposition of the proposition \ref{jdecomp}. In particular $m\leq m'$. So $P(m')$ implies $P(m)$. So we 
can suppose that $\alpha_1(m'),\alpha_2(m'),\alpha_3(m')$ are satisfied.

\noindent If $j\neq i$ : let us prove that $\alpha_1(m),\alpha_2(m),\alpha_3(m)$ are satisfied :

$\alpha_1(m)$ : it follows from the lemma \ref{extl} with $j\neq i$ that we have $L_j(m')=Y_{i,a}L_j(Y_{i,a}^{-1}m')$, 
and so 
$\alpha_1(m')\Rightarrow \alpha_1(m)$.

$\alpha_2(m)$ : as $m(m')^{-1}\in\ZZ[A_{j,b}^{-1}]_{b\in\CC^*}$, we 
have $k(m)=k(m')$ and so $\alpha_2(m')$ implies $\alpha_2(m)$. 

$\alpha_3(m)$ : for $j'\neq j$, we get $\alpha_3(m)$ for $j'$ in the same way. Let us look at it for 
$j$ : for all $l\in\ZZ$, $u_{j,l}(m_{k+1,a}^{(i)})=0$, and so it follows from $\alpha_3(m')$ that for $l< 
2r_i(k(m)+1)$ we have $u_{j,l}(m)=0$. So for such a $l$, $A_{j,aq^{l+r_j}}^{-1}$ does not appear in $L_j(m')$, and we 
get $\alpha_3(m)$.

\noindent If $j=i$ : it follows from $\alpha_3(m')$ that for $0\leq k'\leq k(m')$, $u_{i,aq_i^{2k'}}(m')=1$. So it 
follows from the proposition \ref{aidesldeux} for $L_i(m')$ that for $0\leq k'\leq k(m)$ : 
\begin{equation}\label{argup} w_{i,aq_i^{1+2k'}}(m)\geq 1\Rightarrow \forall k'\leq k''\leq k-1, 
w_{i,aq_i^{1+2k''}}(m)\geq 1.\end{equation} 
So if $k(m)< 1$, $\alpha_2(m')\Rightarrow P(m)$. If $k(m)\geq 1$ the properties 
$\alpha_1(m),\alpha_2(m),\alpha_3(m)$ are satisfied :

$\alpha_1(m)$ : we have $u_{i,a}(m')\geq 1$ and so we can use the (2) of the lemma \ref{extl}. As $v_{i,aq_i}(m(m')^{-1})=0$ we have $\alpha_1(m')\Rightarrow \alpha_1(m)$.

$\alpha_2(m)$ : consequence of the argument (\ref{argup}).

$\alpha_3(m)$ : we have $k(m)\leq k(m')$ and $m(m')^{-1}\in\ZZ[A_{i,b}^{-1}]_{b\in\CC^*}$ and so $\alpha_3(m)$ with 
$j'\neq i$ is clear. For $j'=i$ the property follows from the definition of $k(m)$ and the argument (\ref{argup}).\qed

\subsubsection {End of the proof of the theorem \ref{tsyst} (2) :} Let us denote : 
$$R_{k,a}^{(i)}=(m_{k,a}^{(i)})^{-1}\chi_q(W_{k,a}^{(i)})\in 1+A_{i,aq_i^{2k-1}}^{-1}\ZZ[A_{j,aq^m}^{-1}]_{j\in I, m\geq 0},$$
$$E_{k,a}^{(i)}=(m_{k+1,a}^{(i)})^{-1}A_{i,aq_i}A_{i,aq_i^3}...A_{i,aq_i^{2k+1}}(\chi_q(W_{k+1,a}^{(i)})-Y_{i,a}\chi_q(W_{k,aq_i^2}^{(i)})).$$ 
We have $R_{k+1,a}^{(i)}=R_{k,aq_i^2}^{(i)}+A_{i,aq_i}^{-1}A_{i,aq_i^3}^{-1}...A_{i,aq_i^{2k+1}}^{-1}E_{k,a}^{(i)}$ and from lemma \ref{aideconv} we have $E_{k,a}^{(i)}\in \ZZ[A_{j,aq^m}^{-1}]_{j\in I, m\geq 0}$.
In particular by induction on $k$ we get :
$$R_{k+1,a}^{(i)}=1+\underset{k'=0...k}{\sum}A_{i,aq_i^{1+2(k-k')}}^{-1}A_{i,aq_i^{3+2(k-k')}}^{-1}...A_{i,aq_i^{2k+1}}^{-1}E_{k',aq_i^{2(k-k')}}^{(i)},$$
$$R_{k+1,aq_i^{-2(k+1)}}^{(i)}=1+\underset{k'=0...k}{\sum}A_{i,aq_i^{1+2(-k'-1)}}^{-1}A_{i,aq_i^{3+2(-k'-1)}}^{-1}...A_{i,aq_i^{-1}}^{-1}E_{k',aq_i^{2(-k'-1)}}^{(i)}.$$
We get a graded sum (because all monomial $m$ of the $k'$th term of the sum satisfies $w(m)\geq k'+1$) and so the formal power series $R_{k+1,aq_i^{-2(k+1)}}^{(i)}$ has a limit when $k\rightarrow \infty$.\qed 

\subsubsection{Complement}\label{partdeuxcomp} In this section we give complements that are not used for the main results of this paper : in \cite{hkoty} the asymptotic property is different than the asymptotic property used in the theorem \ref{tsyst} (see \cite{k} for more comments on the property that we use). Although we do not use the property of \cite{hkoty} for the purpose of this paper, we can prove the following version of it :

\begin{prop}\label{compa} Let $i\in I$, let $L$ be the dimension of a fundamental representation $L=\text{dim}(V_{i,a})$ (it is independent of $a\in\CC^*$) and let $r<1$ such that $1/r > 2L$. 

\noindent (1) On the domain $|A_{j,a}|\geq 1/r$ we have :
$\underset{k\rightarrow \infty}{lim}\frac{\chi_q(W_{k,aq_i^{2-2k}}^{(i)})}{\chi_q(W_{k+1,aq_i^{-2k}}^{(i)})}=Y_{i,a}^{-1}$.

\noindent (2) On the domain $|e^{\alpha_j}|\geq 1/r$ we have :
$\underset{k\rightarrow \infty}{lim}\frac{Q_k^{(i)}}{Q_{k+1}^{(i)}}=e^{-\Lambda_i}$.\end{prop}

\noindent The asymptotic property of \cite{hkoty} is the statement (2) on the domain $|e^{\alpha_j}|>1$. But the hypothesis $|e^{\alpha_j}| > 1/r$ is enough for their proof : if a solution of the $Q$-system is a sum of characters of $\Glie$-modules and satisfies the property (2) of the proposition \ref{compa}, then it equals the formulas $\mathcal{F}(\nu)$.

\demo It suffices to prove (1). Consider :
$$\frac{Y_{i,a}\chi_q(W_{k,aq_i^{2-2k}})}{\chi_q(W_{k+1,aq_i^{-2k}})}=\frac{R_{k,aq_i^{-2(k-1)}}^{(i)}}{R_{k+1,aq_i^{-2k}}^{(i)}}
=1-\frac{A_{i,aq_i^{1-2k}}^{-1}A_{i,aq_i^{3-2k}}^{-1}...A_{i,aq_i}^{-1}E_{k,aq_i^{-2k}}^{(i)}}{R_{k+1,aq_i^{-2k}}^{(i)}}$$
$$=1-\frac{A_{i,aq_i^{1-2k}}^{-1}A_{i,aq_i^{3-2k}}^{-1}...A_{i,aq_i}^{-1}E_{k,aq_i^{-2k}}^{(i)}}{1+A_{i,aq_i}^{-1}\underset{k'=0...k}{\sum}A_{i,aq_i^{1-2k'}}^{-1}A_{i,aq_i^{3-2k'}}^{-1}...A_{i,aq_i^{-1}}^{-1}E_{k',aq_i^{-2k'}}^{(i)}}$$
As $E_{k,a}^{(i)}$ has positive coefficients (lemma \ref{aideconv}) we have $|E_{k,a}^{(i)}|\leq \text{dim}(W_{k+1,a}^{(i)}) \leq L^{k+1}$. So :
$$|\frac{Y_{i,a}\chi_q(W_{k,aq_i^{2-2k}})}{\chi_q(W_{k+1,aq_i^{-2k}})}-1|\leq r^{k+1} L^{k+1}\underset{j\geq 0}{\sum}r^j(\underset{k'=0...k}{\sum}r^{k'} L^{k'+1})^j\leq (rL)^{k+1}\underset{j\geq 0}{\sum}(Lr)^j\frac{1}{(1-rL)^j}$$
and the last term has the limit 0 when $k\rightarrow \infty$ because $rL<1/2$ and $rL/(1-rL)<1$.\qed

\noindent Note that we could replace the condition $1/r > 2L$ by the condition $1/r > 2(L-1)$ because the property $P$ in the proof of the lemma \ref{aideconv} implies that $E_{k,a}^{(i)}$ is a sum of monomials of $(\chi_q(W_{1,a}^{(i)})-Y_{i,a})(\chi_q(W_{1,aq_i^2}^{(i)})-Y_{i,aq_i^2})...(\chi_q(W_{1,aq_i^{2k}}^{(i)})-Y_{i,aq_i^{2k}})$, and so $|E_{k,a}^{(i)}|\leq (L-1)^{k+1}$.

\section{Exact sequence}\label{sequ} The category of finite dimensional representations of quantum affine algebras is 
not semi-simple, and so the $T$-system can not be directly written in terms of modules. In \cite{Nad} Nakajima proved 
that the $T$-system can be written in the form of an exact sequence for simply-laced cases. We present here a new proof 
(without $q,t$-characters \cite{Nab}) which allows us to extend the result to non simply-laced cases :

\begin{thm}\label{exact} Let $i\in I, a\in\CC^*, k\geq 1$. We have :

(1) The module $S_{k,a}^{(i)}$ is special and simple.

(2) The module $W_{k+1,a}^{(i)}\otimes W_{k-1,aq_i^2}^{(i)}$ is simple.

(3) There exists an exact sequence : $$0\rightarrow S_{k,a}^{(i)}\rightarrow W_{k,a}^{(i)}\otimes
W_{k,aq_i^2}^{(i)}\rightarrow W_{k+1,a}^{(i)}\otimes W_{k-1,aq_i^2}^{(i)}\rightarrow 0.$$\end{thm}

\noindent The $(1)$ is a direct consequence of the proposition \ref{xspecial}.

\subsection{Proof of the theorem \ref{exact} (2)} Suppose that $W_{k+1,a}^{(i)}\otimes W_{k-1,aq_i^2}^{(i)}$ is not
simple. The dominant monomials of $\chi_q(W_{k+1,a}^{(i)})\chi_q(W_{k-1,aq_i^2}^{(i)})$ are given in the lemma
\ref{premier}. So there is $0\leq R\leq k-2$ such that :
$$\chi_q(W_{k+1,a}^{(i)})\chi_q(W_{k-1,aq_i^2}^{(i)})=\chi_q(L(M))+\chi_q(L(M_R))+\underset{v(m'M^{-1})>R+1}{\sum}n_{m'}\chi_q(L(m'))$$
$$\text{where }M_R=MA_{i,aq_i^{2k-1}}^{-1}...A_{i,aq_i^{2k-1-2R}}^{-1}\text{ and }n_{m'}\geq 0.$$ 
But we have :
$$M_R=Y_{i,a}Y_{i,aq_i^2}^2...Y_{i,aq_i^{2(2k-2(R+2))}}^2Y_{i,aq_i^{2k-2(R+1)}}\underset{b\in\CC^*, j\neq
i}{\prod}Y_{j,b}^{u_{j,b}(M_R)}.$$ 
In particular $M_RA_{i,aq_i^{2k-1-2R}}^{-1}\in\mathcal{M}(L(M_R))$. But the 
monomials of $\mathcal{M}(W_{k+1,a}^{(i)})\mathcal{M}(W_{k-1,aq_i^2}^{(i)})$ lower than $M_RA_{i,b}^{-1}$ ($b\in\CC^*$) are lower than 
$M_RA_{i,aq_i^{2k-3-2R}}^{-1}$ or than $M_RA_{i,aq_i^{2k+1}}^{-1}$, contradiction.\qed

\subsection{Proof of the theorem \ref{exact} (3)} To prove (3) we can adapt arguments of \cite{Nad} :

\begin{thm}\cite{c} Let $i_1,...,i_L\in I$, $a_1,...,a_l\in\CC$ and $m_1,...,m_l\geq 1$. the condition :
$$l<m\Rightarrow \forall p\geq 0, a_l/a_m\neq q^{r_{i_l}k_l-r_{i_m}k_m-r_{i_l}-r_{i_m}-p}$$
implies that $W_{k_1,a_1}^{(i_1)}\otimes ...\otimes W_{k_L,a_L}^{(i_L)}$ is an $l$-highest weight module.\end{thm}

\begin{cor} For $i\in I, a\in\CC^*, k\geq 1$, the module $W_{k,a}^{(i)}\otimes W_{k,aq_i^2}^{(i)}$ is an $l$-highest weight module.\end{cor}

{\it End of the proof of the theorem \ref{exact} (3) :} $W_{k,a}^{(i)}\otimes W_{k,aq_i^2}^{(i)}$ has a unique simple 
quotient. It is isomorphic to the simple $l$-highest weight module $L(m_{k,a}^{(i)}m_{k,aq_i^2}^{(i)})\simeq
W_{k+1,a}^{(i)}\otimes W_{k-1,aq_i^2}^{(i)}$. So it follows from the theorem \ref{tsyst} that the unique maximal proper submodule $\mathcal{M}$ of $W_{k,a}^{(i)}\otimes W_{k,aq_i^2}^{(i)}$ has the $q$-character of $S_{k,a}^{(i)}$. But $S_{k,a}^{(i)}$ is simple, so $\mathcal{M}\simeq S_{k,a}^{(i)}$.\qed

\section{Formulas for the $T$-systems}\label{exptsyst} In this section we give explicit formulas for the $T$-system of 
the theorem \ref{tsyst}. It is written in the form of an exact sequence from theorem \ref{exact} (the $T$-system is 
the $q$-characters identity derived from the exact sequence; for identification 
with the functional formulas of \cite{kns}, let $r=\underset{i\in 
I}{\text{max}}(r_i)$, let identify $u+\ZZ/r$ with $aq^{\ZZ}$ by $u+m/r\rightarrow aq^{m}$, and $T_k^{(i)}(u+m/r)$ with 
$\chi_q(W_{k,aq^{m-r_i(k-1)}}^{(i)})$).

\noindent Type $ADE$ :

$$0\rightarrow 
\underset{j/C_{i,j}=-1}{\bigotimes}W_{k,aq}^{(j)}\rightarrow 
W_{k,a}^{(i)}\otimes W_{k,aq^2}^{(i)}\rightarrow W_{k-1,aq^2}^{(i)}\otimes W_{k+1,a}^{(i)}\rightarrow 0.$$

\noindent Type $B_n$ :

For $2\leq i\leq n-2$ :

$$0\rightarrow W_{k,aq^2}^{(i-1)}\otimes W_{k,aq^2}^{(i+1)}\rightarrow 
W_{k,a}^{(i)}\otimes W_{k,aq^4}^{(i)}\rightarrow W_{k-1,aq^4}^{(i)}\otimes W_{k+1,a}^{(i)}\rightarrow 0,$$

$$0\rightarrow W_{k,aq^2}^{(2)}\rightarrow W_{k,a}^{(1)}\otimes W_{k,aq^4}^{(1)}\rightarrow W_{k-1,aq^4}^{(1)}\otimes 
W_{k+1,a}^{(1)}\rightarrow 0,$$

$$0\rightarrow 
W_{k,aq^2}^{(n-2)}\otimes 
W_{2k,aq}^{(n)}\rightarrow 
W_{k,a}^{(n-1)}\otimes W_{k,aq^4}^{(n-1)}\rightarrow W_{k-1,aq^4}^{(n-1)}\otimes W_{k+1,a}^{(n-1)}\rightarrow 0,$$

$$0\rightarrow W_{r,aq}^{(n-1)}\otimes W_{r,aq^3}^{(n-1)}\rightarrow W_{2r,a}^{(n)}\otimes 
W_{2r,aq^2}^{(n)}\rightarrow W_{2r-1,aq^2}^{(n)}\otimes W_{2r+1,a}^{(n)}\rightarrow 0,$$

$$0\rightarrow 
W_{r+1,aq}^{(n-1)}\otimes W_{r,aq^3}^{(n-1)}\rightarrow W_{2r+1,a}^{(n)}\otimes W_{2r+1,aq^2}^{(n)}\rightarrow 
W_{2r,aq^2}^{(n)}\otimes W_{2r+2,a}^{(n)}\rightarrow 0.$$

\noindent Type $C_n$ :

For $2\leq i\leq n-2$ :

$$0\rightarrow W_{k,aq}^{(i-1)}\otimes W_{k,aq}^{(i+1)}\rightarrow 
W_{k,a}^{(i)}\otimes W_{k,aq^2}^{(i)}\rightarrow W_{k-1,aq^2}^{(i)}\otimes W_{k+1,a}^{(i)}\rightarrow 0,$$

$$0\rightarrow 
W_{k,aq}^{(2)}\rightarrow 
W_{k,a}^{(1)}\otimes W_{k,aq^2}^{(1)}\rightarrow W_{k-1,aq^2}^{(1)}\otimes W_{k+1,a}^{(1)}\rightarrow 0,$$

$$0\rightarrow 
W_{2r,aq}^{(n-2)}\otimes 
W_{r,aq}^{(n)}\otimes 
W_{r,aq^3}^{(n)}\rightarrow 
W_{2r,a}^{(n-1)}\otimes W_{2r,aq^2}^{(n-1)}\rightarrow W_{2r-1,aq^2}^{(n-1)}\otimes W_{2r+1,a}^{(n-1)}\rightarrow 0,$$

$$0\rightarrow W_{2r+1,aq}^{(n-2)}\otimes W_{r+1,aq}^{(n)}\otimes W_{r,aq^3}^{(n)}\rightarrow 
W_{2r+1,a}^{(n-1)}\otimes W_{2r+1,aq^2}^{(n-1)}\rightarrow W_{2r,aq^2}^{(n-1)}\otimes W_{2r+2,a}^{(n-1)}\rightarrow 
0,$$

$$0\rightarrow W_{2k,aq}^{(n-1)}\rightarrow  
W_{k,a}^{(n)}\otimes W_{k,aq^4}^{(n)}\rightarrow W_{k-1,aq^4}^{(n)}\otimes W_{k+1,a}^{(n)}\rightarrow 0.$$

\noindent Type $F_4$ :

$$0\rightarrow 
W_{k,aq^2}^{(2)} 
\rightarrow W_{k,a}^{(1)}\otimes W_{k,aq^4}^{(1)}\rightarrow W_{k-1,aq^4}^{(1)}\otimes 
W_{k+1,a}^{(1)}\rightarrow 0,$$

$$0\rightarrow 
W_{k,aq^2}^{(1)}\otimes 
W_{2k,aq}^{(3)}\rightarrow W_{k,a}^{(2)}\otimes W_{k,aq^4}^{(2)}\rightarrow W_{k-1,aq^4}^{(2)}\otimes 
W_{k+1,a}^{(2)}\rightarrow 
0,$$

$$0\rightarrow 
W_{r,aq}^{(2)}\otimes W_{r,aq^3}^{(2)}\otimes W_{2r,aq}^{(4)}\rightarrow W_{2r,a}^{(3)}\otimes 
W_{2r,aq^2}^{(3)}\rightarrow W_{2r-1,aq^2}^{(3)}\otimes W_{2r+1,a}^{(3)}\rightarrow 0,$$

$$0\rightarrow 
W_{r+1,aq}^{(2)}\otimes 
W_{r,aq^3}^{(2)}\otimes 
W_{2r+1,aq}^{(4)}\rightarrow 
W_{2r+1,a}^{(3)}\otimes W_{2r+1,aq^2}^{(3)}\rightarrow W_{2r,aq^2}^{(3)}\otimes W_{2r+2,a}^{(3)}\rightarrow 0,$$

$$0\rightarrow W_{k,aq}^{(3)}\rightarrow W_{k,a}^{(4)}\otimes W_{k,aq^2}^{(4)}\rightarrow W_{k-1,aq^2}^{(4)}\otimes 
W_{k+1,a}^{(4)}\rightarrow 0.$$

\noindent Type $G_2$ :

$$0\rightarrow W_{3k,aq}^{(2)}\rightarrow W_{k,a}^{(1)}\otimes W_{k,aq^6}^{(1)}\rightarrow W_{k-1,aq^6}^{(1)}\otimes 
W_{k+1,a}^{(1)}\rightarrow 0,$$

$$0\rightarrow W_{r,aq}^{(1)}\otimes W_{r,aq^3}^{(1)}\otimes W_{r,aq^5}^{(1)}\rightarrow 
W_{3r,a}^{(2)}\otimes W_{3r,aq^2}^{(2)}\rightarrow W_{3r-1,aq^2}^{(2)}\otimes W_{3r+1,a}^{(2)}\rightarrow 0,$$

$$0\rightarrow W_{r+1,aq}^{(1)}\otimes W_{r,aq^3}^{(1)}\otimes W_{r,aq^5}^{(1)}\rightarrow
W_{3r+1,a}^{(2)}\otimes W_{3r+1,aq^2}^{(2)}\rightarrow W_{3r,aq^2}^{(2)}\otimes W_{3r+2,a}^{(2)}\rightarrow 0,$$

$$0\rightarrow  
W_{r+1,aq}^{(1)}\otimes 
W_{r+1,aq^3}^{(1)}\otimes 
W_{r,aq^5}^{(1)}\rightarrow 
W_{3r+2,a}^{(2)}\otimes W_{3r+2,aq^2}^{(2)}\rightarrow W_{3r+1,aq^2}^{(2)}\otimes W_{3r+3,a}^{(2)}\rightarrow 0.$$


\begin{thebibliography}{99}

\bibitem[1]{be} {\bf H.A. Bethe}, {Zur Theorie der Metalle. I. Eigenwerte und Eigenfunktionen der linearen Atomkette}

{Z. Phys. 71, 205-226 (1931)}

\mk

\bibitem[2]{bou} {\bf N. Bourbaki}, {\it Groupes et alg\`ebres de Lie} 

{Chapitres IV-VI, Hermann (1968)}

\mk

\bibitem[3]{c0} {\bf V. Chari}, {\it On the fermionic formula and the 
Kirillov-Reshetikhin conjecture} 

{Int. Math. Res. Not. {\bf 2001}, no. 12, 629--654 (2001)}

\mk

\bibitem[4]{c} {\bf V. Chari}, {\it Braid group actions and tensor products}

{Int. Math. Res. Not. {\bf 2002}, no. 7, 357--382 (2002)}

\mk

\bibitem[5]{Cha0} {\bf V. Chari and A. Pressley}, {\it Quantum Affine Algebras}

{Comm. Math. Phys. {\bf 142}, 261-283 (1991)}

\mk

\bibitem[6]{Cha}{\bf V. Chari and A. Pressley}, {\it Quantum affine algebras and their representations}

{dans Representations of groups (Banff, AB, 1994),59-78, CMS Conf. Proc, {\bf 16}, Amer. Math. Soc., Providence, RI (1995)}

\mk

\bibitem[7]{Cha2}{\bf V. Chari and A. Pressley}, {\it A Guide to Quantum Groups}

{Cambridge University Press, Cambridge (1994)} 

\mk

\bibitem[8]{Cha4}{\bf V. Chari and A. Pressley}, {\it Integrable and Weyl modules for quantum affine ${\rm sl}\sb 2$}

{Quantum groups and Lie theory (Durham,
   1999), 48--62, London Math. Soc. Lecture Note Ser., {\bf 290}, Cambridge Univ. Press, Cambridge, (2001)}

\mk

\bibitem[9]{Dri1}{\bf V. G. Drinfel'd}, {\it Quantum groups} 

{Proceedings of the International Congress of Mathematicians, Vol. 1, 2 (Berkeley, Calif., 1986), 798--820, Amer. Math. Soc., Providence, RI, (1987)}

\mk

\bibitem[10]{Fre2} {\bf E. Frenkel and E. Mukhin}, 
{\it Combinatorics of $q$-Characters of Finite-Dimensional Representations of Quantum Affine Algebras} 

{Comm. Math. Phy., vol {\bf 216}, no. 1, pp 23-57 (2001)}

\mk

\bibitem[11]{Fre} {\bf E. Frenkel and N. Reshetikhin}, {\it The $q$-Characters of Representations of Quantum Affine Algebras and Deformations of $W$-Algebras} 

{Recent Developments in Quantum Affine Algebras and related topics, Cont. Math., vol. {\bf 248}, pp 163-205 (1999)}

\mk

\bibitem[12]{her02} {\bf D. Hernandez}, {\it Algebraic approach to q,t-characters}
 
{Adv. Math {\bf 187}, no. 1, 1--52 (2004)}

\mk

\bibitem[13]{her04} {\bf D. Hernandez}, {\it Representations of quantum affinizations and fusion product}

{Transform. Groups {\bf 10}, no. 2, 163--200 (2005)}

\mk

\bibitem[14]{her05} {\bf D. Hernandez}, {\it Monomials of q and q,t-chraracters for non simply-laced quantum 
affinizations}

{Math. Z. {\bf 250}, no. 2, 443--473 (2005)}

\mk

\bibitem[15]{hkott} {\bf G. Hatayama, A. Kuniba, M. Okado, T. Takagi and Z. Tsuboi}, {\it Paths, crystals and fermionic formulae} 

{MathPhys odyssey, 2001, 205--272, Prog. Math. Phys., {\bf 23}, Birkhäuser Boston, Boston, MA, 2002}

\mk

\bibitem[16]{hkoty} {\bf G. Hatayama, A. Kuniba, M. Okado, T. Takagi and Y. Yamada}, {\it Remarks on fermionic formula} 

{in Recent developments in quantum affine algebras
   and related topics (Raleigh, NC, 1998), 243--291, Contemp. Math., {\bf 248}, Amer. Math. Soc., Providence, RI (1999)}

\mk

\bibitem[17]{jim} {\bf M. Jimbo}, {\it A $q$-difference analogue of $\U({\Glie})$ and the Yang-Baxter equation}

{Lett. Math. Phys. {\bf 10}, no. 1, 63--69 (1985)}

\mk

\bibitem[18]{kac} {\bf V. Kac}, {\it Infinite dimensional Lie algebras} 

{3rd Edition, Cambridge University Press (1990)}

\mk

\bibitem[19]{kas} {\bf M. Kashiwara}, {\it On level-zero representations of quantized affine algebras}

{Duke Math. J. {\bf 112}, no. 1, 117--175 (2002)}

\mk

\bibitem[20]{ki1} {\bf A. N. Kirillov}, {\it Combinatorial identities and completeness of states of the Heisenberg magnet}

{Questions in quantum field theory and statistical physics, 4. 
Zap. Nauchn. Sem. Leningrad. Otdel. Mat. Inst. Steklov. (LOMI) {\bf 131} (1983), 88--105 (1983)}

\bibitem[21]{ki2} {\bf A. N. Kirillov}, {\it Completeness of states of the generalized Heisenberg magnet} 

{Automorphic functions and number theory, II. Zap. Nauchn.
   Sem. Leningrad. Otdel. Mat. Inst. Steklov. (LOMI) 134169--189 (1984)}

\mk

\bibitem[22]{ki3} {\bf A. N. Kirillov}, {\it Identities for the Rogers dilogarithmic function connected with simple Lie algebras} 

{J. Soviet Math. {\bf 47}, no. 2, 2450--2459 (1989); translated from Zap. Nauchn. Sem. Leningrad. Otdel. Mat. Inst. Steklov. (LOMI) 164, Differentsialnaya Geom. Gruppy Li i Mekh. IX,  121--133, {\bf 198} (1987)}

\mk

\bibitem[23]{kl} {\bf M. Kleber}, {\it Combinatorial structure of finite-dimensional representations of Yangians: the simply-laced case}

{Internat. Math. Res. Notices {\bf 1997}, no. 4, 187--201 (1997)}

\mk

\bibitem[24]{kn} {\bf H. Knight}, {\it Spectra of tensor products of finite-dimensional representations of Yangians}
   
{J. Algebra {\bf 174}, no. 1, 187--196 (1995)}

\mk

\bibitem[25]{k} {\bf A. Kuniba and T. Nakanishi}, {\it The Bethe equation at $q=0$, the Möbius inversion formula, and
   weight multiplicities. II. The $X_n$ case} 

{J. Algebra {\bf 251}, no. 2, 577--618 (2002)}

\mk

\bibitem[26]{kns} {\bf A. Kuniba, T. Nakanishi and J. Suzuki}, {\it Functional relations in solvable lattice models. I. Functional relations and representation theory}

{Internat. J. Modern Phys. A 9, no. {\bf 30}, 5215--5266 (1994)}

\mk

\bibitem[27]{knt} {\bf A. Kuniba, T. Nakanishi and Z. Tsuboi}, {\it The canonical solutions of the $Q$-systems and the Kirillov-Reshetikhin conjecture} 

{Comm. Math. Phys. {\bf 227}, no. 1, 155--190 (2002)}

\mk

\bibitem[28]{kosy} {\bf A. Kuniba, M. Okado, J. Suzuki, and Y. Yamada}, {\it Difference $L$ operators related to $q$-characters} 

{J. Phys. A {\bf 35}, no. 6, 1415--1435 (2002)}

\mk

\bibitem[29]{kr} {\bf A.N. Kirillov and N. Reshetikhin}, {\it Representations of Yangians and multiplicities of the inclusion of the irreducible components of the tensor product of representations of simple Lie algebras} 

{J. Soviet Math. {\bf 52}, no. 3, 3156--3164 (1990); translated from Zap. Nauchn. Sem. Leningrad. Otdel. Mat. Inst. Steklov. (LOMI) 160, Anal. Teor. Chisel i Teor. Funktsii. 8, 211--221, 301 (1987)}

\mk

\bibitem[30]{Naams} {\bf H. Nakajima}, {\it Quiver varieties and finite-dimensional representations of quantum affine algebras} 

{J. Amer. Math. Soc. {\bf 14}, no. 1, 145--238 (2001)}

\mk

\bibitem[31]{Nab} {\bf H. Nakajima}, {\it Quiver Varieties and $t$-Analogs of $q$-Characters of Quantum Affine Algebras} 

{Ann. of Math. {\bf 160}, 1057--1097 (2004)} 

\mk

\bibitem[32]{Nad}{\bf H. Nakajima}, {\it $t$-analogs of $q$-characters of Kirillov-Reshetikhin modules of quantum 
affine algebras} 

{Represent. Theory {\bf 7}, 259--274 (electronic) (2003)}

\end{thebibliography}
\end{document}